\documentclass{amsart}
\usepackage[alphabetic]{amsrefs}
\usepackage{accents}
\usepackage{graphicx}
\usepackage{amsfonts}
\usepackage{amscd}
\usepackage{amssymb}
\usepackage{xypic}
\xyoption{all}
\usepackage{hyperref}
\usepackage{url}
\usepackage{tikz}
\usepackage{calligra}
\usepackage{extpfeil}
\usetikzlibrary{shapes}
\usetikzlibrary{decorations.markings}
\usepackage{tikz-cd}
\usetikzlibrary{arrows}
\usepackage{multirow}
\usepackage[savepos]{zref}
\newtheorem{thm}{Theorem}[section]
\newtheorem{corollary}[thm]{Corollary}
\newtheorem{lemma}[thm]{Lemma}
\newtheorem{proposition}[thm]{Proposition}

\newtheorem{example}[thm]{Example}

\theoremstyle{definition}
\newtheorem{definition}[thm]{Definition}
\theoremstyle{remark}
\newtheorem{remark}[thm]{Remark}
\numberwithin{equation}{section}
\usepackage[mathscr]{euscript}
\usepackage[T2A,T1]{fontenc}

  {\begin{description}%
    \setlength{\itemsep}{5pt}%
    \setlength{\parskip}{5pt}}%
  {\end{description}}

\newcommand{\form}{\upsilon}

\DeclareMathOperator{\can}{can}

\DeclareMathOperator{\mult}{m}

\DeclareMathOperator{\Pic}{Pic}

\DeclareMathOperator{\SC}{sc}

\DeclareMathOperator{\Comm}{Comm}

\DeclareMathOperator{\Hom}{Hom}

\DeclareMathOperator{\Aut}{Aut}
\DeclareMathOperator{\Int}{Int}

\DeclareMathOperator{\Gal}{Gal}

\DeclareMathOperator{\Ker}{Ker}

\DeclareMathOperator{\val}{val}

\DeclareMathOperator{\Spec}{Spec}

\DeclareMathOperator{\Id}{Id}

\newdir{ >}{{}*!/-5pt/@{>}}
\newdir^{ (}{{}*!/-5pt/@^{(}}

\DeclareMathOperator{\et}{\acute{e}t}

\newcommand{\From}{\colon}
\newcommand{\inar}{\ar@{^{(}->}}
\newcommand{\onar}{\ar@{->>}}

\renewcommand{\Im}{\mathrm{Im}}

\newcommand{\defined}[1]{\underline{{#1}}}
\usepackage{accents}
\newlength{\dtildeheight}

\newcommand{\Baer}{\dotplus}

\newcommand{\Cat}[1]{ {\mathsf{#1}} }
\newcommand{\Fun}[1]{ {\mathsf{#1}} }

\newcommand{\alg}[1]{\boldsymbol{\mathrm{#1}}}
\newcommand{\stack}[1]{{\boldsymbol{\Cat{#1}}}}
\newcommand{\sheaf}[1]{{\mathscr{#1}}}
\newcommand{\shom}{\mathscr{H}\mathit{om}}
\newcommand{\ssym}{\mathscr{S}\mathit{ym}}

\newcommand{\ZZ}{\mathbb Z}

\newcommand{\kk}{\mathfrak{f}}

\renewcommand{\SS}{\mathcal S}

\newcommand{\OO}{\mathcal{O}}

\newcommand{\Into}{\hookrightarrow}
\newcommand{\Onto}{\twoheadrightarrow}
\newcommand{\To}{\rightarrow}

\newcommand{\isom}{\cong}

\newcommand{\jota}{p}

\newcommand{\inarrow}{\arrow[hook]}
\newcommand{\onarrow}{\arrow[two heads]}

\newcommand{\Zero}{{\mathbf 0}}
\makeatletter
\newcommand\@biprod[1]{%
  \vcenter{\hbox{\ooalign{$#1\prod$\cr$#1\coprod$\cr}}}}
\newcommand\biprod{\mathop{\mathpalette\@biprod\relax}\displaylimits}
\makeatother

\newcommand{\defeq}{:=}
\DeclareMathAlphabet{\mathcalligra}{T1}{calligra}{m}{n}
\DeclareMathOperator{\Zar}{Zar}
\DeclareMathOperator{\zar}{zar}

\begin{document}

\title{Covering groups and their integral models}%
\author{Martin H. Weissman}%
\date{\today}

\address{Yale-NUS College, 6 College Ave East, \#B1-01, Singapore 138614}
\email{marty.weissman@yale-nus.edu.sg}%

% \thanks{}

%\keywords{metaplectic, covering, integral model}%
%\date{}%
%\dedicatory{}%
%\commby{}%
% ----------------------------------------------------------------
\begin{abstract}
Given a reductive group $\alg{G}$ over a base scheme $S$, Brylinski and Deligne studied the central extensions of a reductive group $\alg{G}$ by $\alg{K}_2$, viewing both as sheaves of groups on the big Zariski site over $S$.  Their work classified these extensions by three invariants, for $S$ the spectrum of a field.  We expand upon their work to study ``integral models'' of such central extensions, obtaining similar results for $S$ the spectrum of a sufficiently nice ring, e.g., a DVR with finite residue field or a DVR containing a field.  Milder results are obtained for $S$ the spectrum of a Dedekind domain, often conditional on Gersten's conjecture.  

\end{abstract}

\maketitle

\tableofcontents

\section*{Introduction}

\subsection*{Motivation}
Consider a reductive group $\alg{G}$ over a $p$-adic field $F$, with $G = \alg{G}(F)$.  If $G^\circ$ is a hyperspecial maximal compact subgroup of $G$, then one may consider the {\em unramified} representations of $G$ with respect to $G^\circ$.  An unramified irrep of $G$ is an irreducible smooth complex representation of $G$ with a nonzero $G^\circ$-fixed vector.  The classification of unramified irreps by Satake parameters determines the local L-factors of an automorphic representation almost everywhere.  

The choice of hyperspecial maximal compact subgroup can be reformulated algebraically as follows:  an {\em integral model} of $\alg{G}$ will mean a reductive (i.e., smooth, with connected reductive geometric fibres) group over the ring of integers $\OO \subset F$ whose base change to $F$ is endowed with an isomorphism to $\alg{G}$.  It turns out \cite[\S 3.8.1]{Tits} that the $\OO$-points of such a model gives a hyperspecial maximal compact, and all hyperspecial maximal compacts arise as $\OO$-points of such a model.  Therefore, to discuss {\em unramified} representations of a reductive $p$-adic group, it makes sense to begin with a reductive group over $\OO$ instead of $F$.

Now we pass to ``covering groups.''  Just as a ``reductive $p$-adic group'' $G$ arises from the $F$-points of a reductive algebraic group $\alg{G}$ over $F$, a  ``covering'' $\mu_n \Into \tilde G \Onto G$ often arises from an algebraic object:  a central extension of $\alg{G}$ by $\alg{K}_2$ over $F$, as defined by Brylinski and Deligne \cite{B-D}.  We write $\alg{K}_2 \Into \alg{G}' \Onto \alg{G}$ for such a central extension here, and $\mu_n \Into \tilde G \Onto G$ for the resulting central extension of locally compact groups obtained from $F$-points and a Hilbert symbol as in \cite[\S 10.3]{B-D}.

If $G^\circ$ is a hyperspecial maximal compact subgroup of $G$, and $G^\circ \Into \tilde G$ is a splitting of $\tilde G \Onto G$, then one may consider the {\em unramified} (genuine) representations of $\tilde G$ with respect to $G^\circ$.  An unramified irrep of $\tilde G$ is an irreducible smooth complex representation of $\tilde G$ with a nonzero $G^\circ$-fixed vector via this splitting.  A Satake isomorphism exists in this generality by \cite[\S 3.1]{WWL1} (see also \cite[\S 13.10]{McN}).

The hyperspecial maximal compact subgroup $G^\circ \subset G$ may be reformulated algebraically; it arises from an $\OO$-model of $\alg{G}$.  It turns out that splittings $G^\circ \Into \tilde G$ often arise from $\OO$-models of the extension $\alg{K}_2 \Into \alg{G}' \Onto \alg{G}$.  More precisely, an extension $\alg{K}_2 \Into \alg{G}' \Onto \alg{G}$ defined over $\OO$ yields not only a covering group $\mu_n \Into \tilde G \Onto G$, but also a splitting $G^\circ \Into \tilde G$, whenever $n$ is coprime to the residue characteristic of $F$ (i.e., the case of a tame cover).  Therefore, to discuss {\em unramified} representations of a cover of a reductive $p$-adic group, it makes sense to study central extensions of $\alg{G}$ by $\alg{K}_2$ over $\OO$ instead of $F$.  This is the subject of the article.  

\subsection*{Main results}

Let $S$ be a scheme, and let $\alg{G}$ be a reductive group over $S$.  We follow \cite{SGA3} in our conventions, so this means that $\alg{G}$ is a smooth group scheme over $S$ whose geometric fibres are connected reductive groups.  Assume moreover that $\alg{G}$ possesses a maximal torus $\alg{T}$ defined over $S$.  In \cite{B-D}, Brylinski and Deligne study the category central extensions of $\alg{G}$ by $\alg{K}_2$: the category $\Cat{CExt}_S(\alg{G}, \alg{K}_2)$, where $\alg{G}$ and $\alg{K}_2$ are viewed as sheaves of groups on the big Zariski site $S_{\Zar}$.  Such a central extension will be written 
$$\alg{K}_2 \Into \alg{G}' \Onto \alg{G}.$$

When $S = \Spec(F)$, for a field $F$, Brylinski and Deligne classify these central extensions by means of a triple $(Q, \sheaf{D}, f)$ of invariants.  We describe these triples in the Section \ref{BDSection}; they are the objects of a category we call $\Cat{BD}_F(\alg{G},\alg{T})$.  The main result of \cite{B-D} is an equivalence of Picard categories
$$\Fun{BD}_F \From \Cat{CExt}_F(\alg{G}, \alg{K}_2) \xrightarrow{\sim} \Cat{BD}_F(\alg{G}, \alg{T}).$$

Section \ref{ClassSection} reviews this classification, and its much easier $\alg{K}_1$-analogue.  The $\alg{K}_1$-analogue describes the category of central extensions $\alg{G}_{\mult} \Into \alg{G}' \Onto \alg{G}$ of algebraic groups over a field $F$ by means of a category of pairs $(\sheaf{Y}', f)$.  It is phrased as an equivalence of Picard categories
$$\Fun{EZ}_F \From \Cat{CExt}_{F}(\alg{G}, \alg{G}_{\mult}) \xrightarrow{\sim} \Cat{EZ}_F(\alg{G}, \alg{T}).$$

Section \ref{IMSection} is devoted to the theory of {\em integral models} of the central extensions discussed by Brylinski and Deligne.  Suppose that $\OO$ is a Dedekind domain with field of fractions $F$.  Let $S = \Spec(\OO)$, with closed points $S^{(1)}$; for $s \in S^{(1)}$, let $\kk(s)$ denote the corresponding residue field.  Let $\alg{G}$ be a reductive group over $\OO$, and write $\alg{\bar G}_s$ for the special fibre over a closed point $s \in S^{(1)}$.  Let $\eta \From \Spec(F) \To S$ denote the inclusion of the generic point.

There is a Picard category $\Cat{CExt}_\OO(\alg{G}, \alg{K}_2)$ of central extensions of $\alg{G}$ by $\alg{K}_2$ defined over $\OO$.    The construction from \cite[\S 12.11]{B-D} gives a functor 
$$\partial_s \From\Cat{CExt}_F(\alg{G}_F, \alg{K}_2) \To \Cat{CExt}_{\kk(s)}(\alg{\bar G}_s, \alg{\bar G}_{\mult}),$$ 
for all $s \in S^{(1)}$.  The main result of Section \ref{IMSection} is a left-exact sequence of Picard categories,
$$0 \To \Cat{CExt}_\OO(\alg{G}, \alg{K}_2) \xrightarrow{\eta^\ast} \Cat{CExt}_F(\alg{G}_F, \alg{K}_2) \xrightarrow{\partial} \bigoplus_{s \in S^{(1)}} \Cat{CExt}_{\kk(s)}(\alg{\bar G}_s, \alg{\bar G}_{\mult}).$$
In particular, given a central extension $\alg{G}_F' \in \Cat{CExt}_F(\alg{G}_F, \alg{K}_2)$, giving an $\OO$-model of $\alg{G}_F'$ is the same as giving trivializations of the {\em residual extensions} $\partial_s \alg{G}_F'$ for all $s \in S^{(1)}$.  This follows and generalizes \cite[Remark 12.14(iii)]{B-D}.

This result relies on Gersten's conjecture (in weight two) for smooth schemes of finite type over $\OO$; when $\OO$ is a discrete valuation ring (DVR) with finite residue field, a DVR containing a field, or $\OO$ is the ring of $\SS$-integers in a global field of prime characteristic, the necessary cases of Gersten's conjecture hold.  The results of Section \ref{IMSection} also provide an $\OO$-integral version of the functor $\Fun{BD}_F$,
$$\Fun{BD}_\OO \From \Cat{CExt}_\OO(\alg{G}, \alg{K}_2) \To \Cat{BD}_{\OO}(\alg{G}, \alg{T}).$$

Section \ref{RESection} assembles the following diagram of Picard categories and additive functors, with exact rows, and equivalences along vertical arrows.  This section also provides natural isomorphisms making this diagram commute in the 2-categorical sense.
$$\begin{tikzcd}
0 \arrow{r} & \Cat{CExt}_{\OO}(\alg{G}, \alg{K}_2) \arrow{r} \arrow{d}{\Fun{BD}_\OO}[swap]{\sim} & \Cat{CExt}_F(\alg{G}_F, \alg{K}_2) \arrow{r}{\partial} \arrow{d}{\Fun{BD}_F}[swap]{\sim} & \Cat{CExt}_{\kk}(\alg{\bar G}, \alg{\bar G}_{\mult}) \arrow{d}{\Fun{EZ}_\kk}[swap]{\sim} \\
0 \arrow{r} & \Cat{BD}_{\OO}(\alg{G}, \alg{T}) \arrow{r} & \Cat{BD}_F(\alg{G}_F, \alg{T}_F) \arrow{r}{\Fun{val}} & \Cat{EZ}(\alg{\bar G}, \alg{\bar T})
\end{tikzcd}$$
The functor $\Fun{val}$ (along with the natural isomorphism expressing commutativity of a square) provides an effective description of the functor $\partial$ of \cite[\S 12.11]{B-D}.

The equivalence $\Fun{BD}_\OO$ provides a classification of central extensions of $\alg{G}$ by $\alg{K}_2$, over the ring $\OO$.  This extends the main results of \cite{B-D}, at least to the case of DVRs with finite residue field, or DVRs which contain a field.

\section*{Acknowledgments}

This article is part of a larger project on ``covering the Langlands program,'' which was greatly assisted by a workshop at the American Institute of Mathematics.  We thank AIM for their support, for bringing together so many colleagues to discuss covering groups.

In particular, we thank Gordan Savin for his endless advice on covering groups, and motivating many of the results here with his study of unramified representations of covering groups.  We thank Wee Teck Gan and 
Gao Fan for many discussions throughout the writing of this paper.

\section*{Preliminaries}

\subsection*{Sheaves}

Let $S$ be a Noetherian (to play it safe) scheme.  We write $S_{\et}$ for the \'etale site and $S_{\zar}$ for the Zariski site.  Write $S_{\Zar}$ for the {\em big} Zariski site of schemes of finite type over $S$.  We work frequently in the topoi of sheaves on these sites:  $\Cat{Sh}(S_{\et})$ and $\Cat{Sh}(S_{\zar})$ and $\Cat{Sh}(S_{\Zar})$.

We use a cursive font, as in $\sheaf{F}$, for a sheaf on $S_{\et}$.  When $U \To S$ is an \'etale morphism, we write $\sheaf{F}[U]$ for the sections of $\sheaf{F}$ over $U$.  We write $\sheaf{G}_{\mult}$ for the multiplicative group viewed as a sheaf on $S_{\et}$.  

We use a boldface font, as in $\alg{F}$, for a sheaf on the big Zariski site $S_{\Zar}$.  Thus we write $\alg{G}_{\mult}$ (or $\alg{G}_{\mult/S}$) for the multiplicative group, viewed as such a sheaf.  We identify schemes of finite type over $S$ with the Zariski sheaves they represent.  If $\alg{F}$ is a sheaf on $S_{\Zar}$, and $\alg{X}$ is a scheme of finite type over $S$, then we obtain a sheaf $\alg{F}_{\alg{X}}$ on the Zariski site $\alg{X}_{\zar}$.  Often we omit the subscript, and simply view $\alg{F}$ as a system, varying functorially with $\alg{X} \To S$, of Zariski sheaves.  

$S$ will almost always denote the spectrum of a Dedekind domain.  In this case, the \defined{Picard group} of $S$ is the group of equivalence classes of line bundles on $S$, identified with the cohomology groups below.
$$\Pic(S) = H_{\et}^1(S, \sheaf{G}_{\mult}) = H_{\zar}^1(S, \alg{G}_{\mult}).$$

\subsection*{Group schemes}

We follow \cite{SGA3} in our conventions for group schemes over $S$.  In particular, a \defined{reductive group} over $S$ will mean a smooth group scheme $\alg{G} \To S$ whose geometric fibres are connected reductive algebraic groups.  Similarly, a \defined{torus} $\alg{T} \To S$ will mean a smooth group scheme, whose geometric fibres are algebraic tori.  

When $\alg{T} \To S$ is a torus, we view its characters $\sheaf{X} = \shom(\alg{G}_{\mult}, \alg{T})$ and cocharacters $\sheaf{Y} = \shom(\alg{G}_{\mult}, \alg{T})$ as local systems on $S_{\et}$.  When $\alg{T} \subset \alg{G}$ is a maximal torus in a reductive group over $S$, the Weyl group will be viewed as a sheaf $\sheaf{W}$ of finite groups on $S_{\et}$.  

When $\alg{T} \To \Spec(\OO)$ is a torus over a ring $\OO$, and $y \in \sheaf{Y}[\OO]$, we obtain a group homomorphism $y \From \OO^\times \To \alg{T}(\OO)$.  We write this ``exponentially'' as follows:
$$y \From \OO^\times \To \alg{T}(\OO), \quad u \mapsto u^y.$$

\subsection*{K-theory of schemes}

We refer to Quillen \cite{Quillen} and Bloch \cite[Chapter 4]{BlochAlgCycle}, for key facts in the K-theory of rings and schemes.  When $\alg{X} \To S$ is a scheme of finite type, and $\alg{U} \subset \alg{X}$ an affine open, $\alg{U} = \Spec(A)$, Quillen's algebraic K-theory provides abelian groups $K_i(\alg{U}) \defeq K_i(A)$ for $i \geq 0$.  We write $\alg{K}_i$ for the Zariski sheaf on $\alg{X}$ associated to the presheaf $\alg{U} \mapsto K_i(\alg{U})$.  As $\alg{X}$ varies over schemes of finite type over $S$, these $\alg{K}_i$ form sheaves of abelian groups on $S_{\Zar}$.  We will only use $\alg{K}_0$, $\alg{K}_1$, and $\alg{K}_2$ in this paper.

Note that $\alg{K}_0$ is the constant sheaf $\ZZ$, and $\alg{K}_1 = \alg{G}_{\mult}$, since $K_1(A) = A^\times$ for any local ring $A$.  For any local ring $A$, the identification $K_1(A) = A^\times$ gives a $\ZZ$-bilinear pairing
$$\{ \bullet, \bullet \} \From A^\times \times A^\times \To K_2(A).$$
This is called the Steinberg symbol, and it satisfies $\{ a,-a \} = 1$, $\{ a, 1-a \} = 1$, and $\{ a_1, a_2 \} \{ a_2, a_1 \} = 1$, whenever $a, a_1, a_2 \in A^\times$.  When $F$ is a field, these relations suffice to characterize $K_2(F)$ as a quotient of $F^\times \otimes_\ZZ F^\times$.

When $\OO$ is a Dedekind domain and $S = \Spec(\OO)$, we write $F$ for its fraction field, and $\eta \From \Spec(F) \Into \Spec(\OO)$ for the generic point of $S$.  Write $S^{(1)}$ for the set of points of $S$ of codimension $1$, i.e., the set of maximal ideals.  For any $s \in S^{(1)}$, write $\kk(s)$ for the associated residue field.  In this level of generality, Quillen \cite[Corollary, p.113]{Quillen} finds a long exact sequence of K-theory groups,
$$K_{i+1}(F) \To \bigoplus_{s \in S^{(1)}} K_i( \kk(s)) \To K_i(\OO) \To K_i(F) \To \bigoplus_{s \in S^{(1)}} K_{i-1}(\kk(s)) \To \cdots.$$

Two special cases of this long exact sequence will arise repeatedly.  First, when $\OO$ is a DVR with residue field $\kk$, we find a short exact sequence
\begin{equation}
\label{K1K0SequenceLocal}
K_1(\OO) \Into K_1(F) \Onto K_0(\kk).
\end{equation}
In this setting $K_1(\OO) = \OO^\times$, $K_1(F) = F^\times$, $K_0(\kk) = \ZZ$, and the map $K_1(F) \To K_0(\kk)$ is the valuation (normalized always to have $\val(F^\times) = \ZZ$).

Second, consider the ring $\OO_\SS$ of $\SS$-integers in a global field $F$.  Writing $S = \Spec(\OO_\SS)$, we find that $S^{(1)}$ is the set of maximal ideals of $\OO_\SS$ outside of $\SS$.  Then every residue field $\kk(s)$ of $\OO_\SS$ is finite, and so $K_2(\kk(s)) = 0$.  The main result of Bass-Milnor-Serre \cite[Theorem 4.1]{BMS} implies that $K_1(\OO_\SS) = \OO_\SS^\times$, from which it follows that $K_1(\OO_\SS) \To K_1(F)$ is injective.  This gives a short exact sequence
\begin{equation}
\label{K2K1SequenceGlobal}
K_2(\OO_\SS) \Into K_2(F) \Onto \bigoplus_{s \in S^{(1)}} \kk(s)^\times.
\end{equation}

When the set of places $\SS$ is sufficiently large, the Picard group $\Pic(\OO_\SS)$ is trivial, from which it follows that $K_0(\OO_\SS) = \ZZ$.  Thus for $\SS$ sufficiently large, we find a global counterpart to the sequence \eqref{K1K0SequenceLocal}
\begin{equation}
\label{K1K0SequenceGlobal}
K_1(\OO_\SS) \Into K_1(F) \Onto \bigoplus_{s \in S^{(1)}} K_0(\kk(s)).
\end{equation}

\subsection*{Central extensions in topoi}

Let $\Cat{T}$ be the topos of sheaves on a site.  When $A$ is an abelian group in $\Cat{T}$, we write $\Cat{Tors}(A)$ for the category of $A$-torsors in $\Cat{T}$.  For two such torsors $R_1, R_2$, their contraction is denoted
$$R_1 \wedge_A R_2 = \frac{ R_1 \times R_2}{(a \ast r_1, r_2) = (r_1, a \ast r_2)}.$$

We refer to \cite{SGA7} and \cite{B-D} for a complete treatment of central extensions of groups in a (Grothendieck) topos.  For any topos $\Cat{T}$, and groups $G$, $A$ in the topos with $A$ abelian, there is a category whose objects are central extensions
$$A \Into E \Onto G$$
of groups in $\Cat{T}$.  Such extensions may be viewed as multiplicative $A_G$-torsors, as in Breen \cite{Bre}.  We write $\Cat{CExt}(G, A)$ for the category of central extensions of $G$ by $A$.  For $E_1, E_2 \in \Cat{CExt}(G,A)$, we write $E_1 \Baer E_2$ for their Baer sum.  This is the central extension, whose associated $A_G$-torsor structure is the contraction of the torsors for $E_1$ and $E_2$.

Given $E \in \Cat{CExt}(G,A)$, lifting followed by conjugation yields an action $\Int \From G \times E \To E$.  If $G$ is abelian too, the commutator provides an alternating form $\Comm \From G \times G \To A$.  The category $\Cat{Ext}(G,A)$ is the full subcategory of $\Cat{CExt}(G, A)$ consisting of {\em abelian} extensions of $G$ by $A$ in the topos $\Cat{T}$, i.e., those in which $\Comm$ is trivial.

\subsection*{Picard categories}

We follow Deligne \cite[Exp. XVIII, \S 1.4]{SGA4T3} in our definitions and treatment of Picard categories, additive functors, and natural transformations.  Here we always assume our Picard categories to be strictly commutative.  Such categories have been studied extensively, sometimes using the term ``symmetric 2-group'' at times.  Some homological algebra, replacing abelian groups by Picard categories, has been developed by C. Bertolin \cite{Ber}, D. Bourn and E.M. Vitale \cite{B-V}, and K.-H. Ulbrich \cite{Ulb}, among others.  

When we write ``$\Cat{P}$ is a Picard category,'' we implicitly mean that a category $\Cat{P}$ endowed with a monoidal functor $\Baer_{\Cat{P}}$, and natural isomorphisms $\Fun{comm}_{\Cat{P}}$ and $\Fun{ass}_{\Cat{P}}$, is a strictly commutative Picard category.

When $A$ is an abelian group in a topos $\Cat{T}$ as before, the category $\Cat{Tors}(A)$ of $A$-torsors forms a Picard category, with monoidal structure given by contraction.  More generally (or by transport of structure), when $G$ is a group in $\Cat{T}$, the category $\Cat{CExt}(G,A)$ of central extensions of $G$ by $A$ forms a Picard category with monoidal structure given by the Baer sum.  

Let $\Cat{Pic}$ be the 2-category, whose objects are Picard categories, where for any such objects $\Cat{X},\Cat{Y}$, the category $\Cat{Hom}(\Cat{X},\Cat{Y})$ consists of additive functors from $X$ to $Y$ and natural transformations among them.  Additive functors may be added; in this way $\Cat{Hom}(\Cat{X},\Cat{Y})$ is a Picard category whenever $\Cat{X}$ and $\Cat{Y}$ are Picard categories.

Following \cite[\S 3]{Ber}, a sequence of Picard categories and additive functors
$$0 \To \Cat{P} \xrightarrow{\alpha} \Cat{Q} \xrightarrow{\beta} \Cat{R}.$$
is called exact if the following conditions hold:
\begin{enumerate}
\item
The composition $\beta \circ \alpha$ is naturally isomorphic to the zero functor;
\item
The functor $\alpha$ induces an equivalence of Picard categories from $\Cat{P}$ to the category $\Cat{Ker}(\beta)$, the category of pairs $(q,s)$ where $q$ is an object of $\Cat{Q}$ and $s$ is an isomorphism from $\beta q$ to $\Zero$.
\end{enumerate}

This notion of exactness is relevant to this paper through the following result.
\begin{proposition}
\label{LExact}
Suppose that $A_1 \xhookrightarrow{\alpha} A_2 \xtwoheadrightarrow{\beta} A_3$ is an exact sequence of abelian groups (in the topos $\Cat{T}$).  Then contraction yields a left-exact sequence of Picard categories and additive functors
$$0 \To \Cat{Tors}(A_1) \To \Cat{Tors}(A_2) \To \Cat{Tors}(A_3).$$
If $G$ is a group in $\Cat{T}$, then this defines a left-exact sequence of Picard categories and additive functors
$$0 \To \Cat{CExt}(G, A_1) \To \Cat{CExt}(G, A_2) \To \Cat{CExt}(G, A_3).$$
\end{proposition}
\proof
Given the exact sequence $A_1 \xhookrightarrow{\alpha} A_2 \xtwoheadrightarrow{\beta} A_3$, we have additive functors of Picard categories,
$$\Cat{Tors}(A_1) \xrightarrow{\alpha} \Cat{Tors}(A_2) \xrightarrow{\beta} \Cat{Tors}(A_3),$$
given by contraction of torsors.  For example, if $R$ is an $A_1$-torsor, then 
$$\alpha(R) = R \wedge_{A_1} A_2 \defeq \frac{ R \times A_2 }{(a_1 \ast r, a_2) \sim (r, \alpha(a_1) \cdot a_2)}.$$

If $R$ is an $A_1$-torsor, then $(R \wedge_{A_1} A_2) \wedge_{A_2} A_3$ is naturally isomorphic to $R \wedge_{A_1} A_3$, where the homomorphism $\beta \alpha \From A_1 \To A_3$ is the trivial homomorphism.  But $R \wedge_{A_1} A_3$ is equal to $(R / A_1) \times A_3$, which is naturally isomorphic to the trivial $A_3$-torsor.  This provides a natural isomorphism $\beta \alpha \xRightarrow{\sim} \Zero$ to the zero functor of Picard categories.

It furthermore provides a functor from $\Cat{Tors}(A_1)$ to the category of pairs $(R_2, \tau)$ where $R_2$ is an $A_2$-torsor and $\tau$ is a neutralization of $\beta R_2$.  To check that this is an equivalence, begin with such a pair $(R_2, \tau)$.  Thus $R_2$ is an $A_2$-torsor, and $\tau \From R_2 \wedge_{A_2} A_3 \To A_3$ is a neutralization of the $A_3$-torsor $\beta R_2$.  

The preimage of $1 \in A_3$ via
$$R_2 \xrightarrow{r \mapsto (r,1)} R_2 \times A_3 \To R_2 \wedge_{A_2} A_3 \xrightarrow{\tau} A_3,$$
is an $A_1$-torsor $R_1 \subset R_2$.  The reader may check that $\alpha R_1$ is then  naturally isomorphic to $R_2$ again.  In this way $\alpha$ provides an equivalence of Picard categories from $\Cat{Tors}(A_1)$ to the category of pairs $(R_2, \tau)$ mentioned above.  

Viewing central extensions as multiplicative bitorsors on $G$, the rest of the proposition follows by transport of structure.
\qed

\section{Brylinski-Deligne invariants}
\label{BDSection}

In this section, fix a field $F$, a reductive group $\alg{G}$ over $F$ and a maximal $F$-torus $\alg{T} \subset \alg{G}$.  Also fix a central extension
$$\alg{K}_2 \Into \alg{G}' \Onto \alg{G},$$
of groups on $F_{\Zar}$.  If $L$ is a field containing $F$, then the vanishing of $H_{\zar}^1(L, \alg{K}_2)$ implies that taking global sections over $L$ gives a central extension of groups,
$$\alg{K}_2(L) \Into \alg{G}'(L) \Onto \alg{G}(L).$$
Recall that the character and cocharacter lattices of $\alg{T}$, viewed as local systems of abelian groups on $F_{\et}$, are written $\sheaf{X}$ and $\sheaf{Y}$.

\subsection{The first Brylinski-Deligne invariant}

To the central extension $\alg{K}_2 \Into \alg{G}' \Onto \alg{G}$, Brylinski and Deligne \cite[\S 3.9]{B-D} associate a Weyl-invariant quadratic form $Q \From \sheaf{Y} \To \ZZ$.  In other words, 
$$Q \in H_{\et}^0(F, \ssym^2 (\sheaf{X})^{\sheaf{W}}).$$

Write $B_Q \From \sheaf{Y} \otimes \sheaf{Y} \To \ZZ$ for the bilinear form associated to $Q$, $B_Q(y_1, y_2) = Q(y_1+y_2) - Q(y_1) - Q(y_2)$.  It satisfies the following, from \cite[Proposition 3.13]{B-D}.
\begin{proposition}
Let $L$ be \textbf{any field} containing $F$ over which $\alg{T}$ splits.  Let $\alg{T}'(L)$ be the resulting central extension,
$$\alg{K}_2(L) \Into \alg{T}'(L) \Onto \alg{T}(L).$$
Then the commutator of this extension satisfies
$$\Comm \left( u_1^{y_1}, u_2^{y_2} \right) = \{ u_1, u_2 \}^{B_Q(y_1, y_2)},$$
for all $y_1, y_2 \in \Hom_L(\alg{G}_{\mult}, \alg{T})$ and all $u_1, u_2 \in L^\times$.  
\end{proposition}

We call the quadratic form $Q$ the \defined{first invariant} of the object $\alg{G}' \in \Cat{CExt}_F(\alg{G}, \alg{K}_2)$.  The first invariant controls the commutators for extensions of tori by $\alg{K}_2$.  When $\alg{G}$ is a simply-connected semisimple group, \cite[Theorems 4.7, 7.2]{B-D} implies that the first invariant classifies central extensions of $\alg{G}$ by $\alg{K}_2$.  
\begin{thm}
\label{UniqueCoversSCSS}
Suppose that $\alg{G}$ is simply-connected, semisimple, over a field $F$, with maximal $F$-torus $\alg{T}$.  Then the central extensions $\alg{K}_2 \Into \alg{G}' \Onto \alg{G}$ have no automorphisms except the identity.  They are classified, up to unique isomorphism, by $\sheaf{W}$-invariant quadratic forms $Q \From \sheaf{Y} \To \ZZ$.  
\end{thm}

From the theorem comes a definition.
\begin{definition}
Given a simply-connected semisimple group $\alg{G}$ over $F$ with maximal $F$-torus $\alg{T}$, and Weyl-invariant quadratic form $Q \From \sheaf{Y} \To \ZZ$, write $\alg{G}_Q'$ for the unique (up to unique isomorphism) central extension of $\alg{G}$ by $\alg{K}_2$ with first invariant $Q$.  It will be called the \defined{canonical} Brylinski-Deligne extension associated to $Q$.
\end{definition}

\begin{remark}
When $\alg{G}$ is simple, split, and $Q$ takes the value $1$ on short coroots, \cite[Proposition 4.15]{B-D} demonstrates that $\alg{G}_Q'$ coincides with Matsumoto's universal central extension \cite{Matsumoto}.\end{remark}

\begin{example}
\label{Ec1}
Let $\alg{G} = \alg{G}_{\mult}$.  Then, for every integer $c$, one may construct a central extension $\alg{K}_2 \Into \alg{E}_c \Onto \alg{G}_{\mult}$ in $\Cat{CExt}_F(\alg{G}_{\mult}, \alg{K}_2)$ as follows:  as a $\alg{K}_2$-torsor on $\alg{G}_{\mult}$, it is trivial.  The multiplicative structure is given by the cocycle $\alg{G}_{\mult} \times \alg{G}_{\mult} \To \alg{K}_2$,
$$u,v \mapsto \{ u,v \}^c.$$
We call this the extension of $\alg{G}_{\mult}$ \defined{incarnated by} $c$.

In particular, the $F$-points of this extension are given by $\alg{E}_c(F) = F^\times \times \alg{K}_2(F)$, with group law,
$$ (u, \alpha) \cdot (v, \beta) = (uv, \alpha \beta \cdot \{ u,v \}^c ).$$
The first invariant of $\alg{E}_c$ is the quadratic form $Q \From \ZZ \To \ZZ$ which satisfies $Q(1) = c$.
\end{example}

\subsection{The second invariant}

Next, associated to $\alg{K}_2 \Into \alg{G}' \Onto \alg{G}$, an object of $\Cat{CExt}_F(\alg{G}, \alg{K}_2)$, and a maximal $F$-torus $\alg{T} \subset \alg{G}$, Brylinski and Deligne construct an extension of sheaves of groups on $F_{\et}$,
$$\sheaf{G}_{\mult} \Into \sheaf{D} \Onto \sheaf{Y}.$$  
Their construction in \cite[\S 3.10]{B-D} proceeds in the followng steps, beginning with a finite separable splitting field $L/F$ of $\alg{T}$.  We work here with Zariski sheaves on $\alg{G}_{\mult/L} = \Spec(L[\form^{\pm 1}])$.
\begin{enumerate}
\item
They observe (following \cite{She}) that $H_{\zar}^1(\alg{G}_{\mult/L} , \alg{K}_2) = 0$, and so taking global sections yields a short exact sequence
$$H_{\zar}^0 \left( \alg{G}_{\mult/L} , \alg{K}_2 \right) \Into H_{\zar}^0 \left(\alg{G}_{\mult/L} , \alg{T}' \right)  \Onto H_{\zar}^0 \left( \alg{G}_{\mult/L} ,\alg{T} \right).$$
\item
Applying another result of \cite{She}, we have $H_{\zar}^0 \left( \alg{G}_{\mult/L} , \alg{K}_2 \right) = K_2(L) \oplus K_1(L)$, giving a residue map $\partial \From H_{\zar}^0 \left( \alg{G}_{\mult/L} , \alg{K}_2 \right) \To L^\times$.  
Pushing out via $\partial$ yields
$$L^\times \Into \partial_\ast  H_{\zar}^0 \left(\alg{G}_{\mult/L} , \alg{T}' \right) \Onto H_{\zar}^0 \left( \alg{G}_{\mult/L} ,\alg{T} \right).$$
\item
There is a canonical injective homomorphism $h \From \sheaf{Y}[L] = \Hom_L(\alg{G}_{\mult}, \alg{T}) \Into H_{\zar}^0 \left( \alg{G}_{\mult/L} ,\alg{T} \right)$.  Pulling back defines a group $\sheaf{D}[L] \defeq h^\ast \partial_\ast H_{\zar}^0 \left(\alg{G}_{\mult/L} , \alg{T}' \right)$, fitting into an extension
$$L^\times \Into \sheaf{D}[L] \Onto \sheaf{Y}[L].$$
(In \cite{B-D}, the group $\sheaf{D}[L]$ is called $\mathcal{E}$ instead.)
\end{enumerate}
Functoriality of this construction, with respect to morphisms $L_1 \To L_2$ of splitting fields of $\alg{T}$, yields a central extension of sheaves of groups on $F_{\et}$,
$$\sheaf{G}_{\mult} \Into \sheaf{D} \Onto \sheaf{Y}.$$
We call this extension the \defined{second invariant} of $\alg{G}'$.  It depends on the maximal $F$-torus $\alg{T}$, but in a predictable way according to \cite[\S 11.12]{B-D}.  

According to \cite[Proposition 3.11]{B-D}), the commutator of this extension is the alternating bilinear map $\Comm \From \sheaf{Y} \otimes \sheaf{Y} \To \sheaf{G}_{\mult}$, given by $\Comm(y_1, y_2) = (-1)^{B_Q(y_1, y_2)}$, with $Q$ the first invariant.    

The second invariant defines a functor
$$\Cat{CExt}_F(\alg{G}, \alg{K}_2) \To \Cat{CExt}_{F_{\et}}(\sheaf{Y}, \sheaf{G}_{\mult}).$$

\begin{definition}
Given a simply-connected semisimple group $\alg{G}$ over $F$ with maximal $F$-torus $\alg{T}$, and Weyl-invariant quadratic form $Q \From \sheaf{Y} \To \ZZ$, we have defined the canonical Brylinski-Deligne extension $\alg{G}_Q'$ by Theorem \ref{UniqueCoversSCSS}.  Define $\sheaf{D}_Q$ to be the second invariant of $\alg{G}_Q'$; it is a central extension
$$\sheaf{G}_{\mult} \Into \sheaf{D}_Q \Onto \sheaf{Y}.$$
The extension $\sheaf{D}_Q$ is characterized by other means in \cite[\S 11]{B-D}.
\end{definition}

\begin{example}
\label{Ec2}
Let $c$ be an integer and let $\alg{K}_2 \Into \alg{E}_c \Onto \alg{G}_{\mult}$ be the extension incarnated by $c$ in Example \ref{Ec1}.  The second invariant is a central extension
$$\sheaf{G}_{\mult} \Into \sheaf{D}_c \Onto \ZZ.$$
Tracing through the construction, we find that $\sheaf{D}_c = \ZZ \times \sheaf{G}_{\mult}$ as a sheaf of sets on $F_{\et}$, with (abelian) group law given by
$$(m, u) \cdot (n,v) = (m + n, uv \cdot (-1)^{mn} ).$$
\end{example}

\subsection{The third invariant}

To define the third invariant, let $\jota \From \alg{G}_{\SC} \To \alg{G}$ denote the simply-connected cover of the derived subgroup.  Pulling back the extension $\alg{K}_2 \Into \alg{G}' \Onto \alg{G}$, and the maximal torus $\alg{T}$, via $\alg{G}_{\SC} \To \alg{G}$ yields an extension $\alg{K}_2 \Into \alg{G}_{\SC}' \Onto \alg{G}_{\SC}$ in $\Cat{CExt}_F(\alg{G}_{\SC}, \alg{K}_2)$ and a maximal torus $\alg{T}_{\SC}$ in $\alg{G}_{\SC}$.  Write $\sheaf{Y}_{\SC}$ for its cocharacter lattice.  The induced map $\jota \From \sheaf{Y}_{\SC} \Into \sheaf{Y}$ is the inclusion of the coroot lattice.    

The construction of the second invariant gives a commutative diagram of sheaves of groups on $F_{\et}$.
$$\begin{tikzcd}
\sheaf{G}_{\mult} \inarrow{r} \arrow{d}{=} &  \sheaf{D}_{\SC} \onarrow{r} \inarrow{d}{\jota} & \sheaf{Y}_{\SC} \inarrow{d}{\jota} \\
\sheaf{G}_{\mult} \inarrow{r} & \sheaf{D} \onarrow{r} & \sheaf{Y}
\end{tikzcd}$$

Theorem \ref{UniqueCoversSCSS} defines a unique isomorphism $f \From \alg{G}_Q' \To \alg{G}_{\SC}'$ in $\Cat{CExt}_F(\alg{G}_{\SC}, \alg{K}_2)$, where $\alg{G}_Q'$ is the canonical extension with first invariant $Q$.  Functoriality of the second invariant provides an isomorphism $f_{\SC}$ in $\Cat{CExt}_{F_{\et}}(\sheaf{Y}_{\SC}, \sheaf{G}_{\mult})$.
$$\begin{tikzcd}
\sheaf{G}_{\mult} \inarrow{r} \arrow{d}{=} &  \sheaf{D}_Q \onarrow{r} \arrow{d}{f_{\SC}} & \sheaf{Y}_{\SC} \arrow{d}{=} \\
\sheaf{G}_{\mult} \inarrow{r} & \sheaf{D}_{\SC} \onarrow{r} & \sheaf{Y}_{\SC}
\end{tikzcd}$$

The \defined{third invariant} of Brylinski and Deligne arises from assembling these two commutative diagrams; the third invariant is the homomorphism $f = \jota \circ f_{\SC}$ of sheaves of groups on $F_{\et}$ which fits into the commutative diagram below.
\begin{equation}
\label{ThirdInvariantDiagram}
\begin{tikzcd}
\sheaf{G}_{\mult} \inarrow{r} \arrow{d}{=} &  \sheaf{D}_Q \onarrow{r} \inarrow{d}{f} & \sheaf{Y}_{\SC} \inarrow{d}{\jota} \\
\sheaf{G}_{\mult} \inarrow{r} & \sheaf{D} \onarrow{r} & \sheaf{Y}
\end{tikzcd}
\end{equation}

\section{Classifications}
\label{ClassSection}
Let $\alg{G}$ be a reductive group over a field $F$ with maximal $F$-torus $\alg{T}$.  The three invariants of the previous section suffice for the classification of central extensions of $\alg{G}$ by $\alg{K}_2$.  We review this result here, using the language of Picard stacks.

\subsection{Brylinski-Deligne classification}

Recall that $\Cat{CExt}_F(\alg{G}, \alg{K}_2)$ is the category of central extensions of $\alg{G}$ by $\alg{K}_2$, viewing both as sheaves of groups on $F_{\Zar}$.  For any finite separable extension $L/F$, we have the corresponding category $\Cat{CExt}_L(\alg{G}_L, \alg{K}_2)$.  This system of categories forms a \defined{Picard stack} (\cite[Expose XVIII, \S 1.4]{SGA4T3} on $F_{\et}$, which we call $\stack{CExt}(\alg{G}, \alg{K}_2)$.  Indeed, the Baer sum provides the monoidal structure and \cite[Theorem 2.7]{B-D} verifies the effectiveness of Galois descent.  
\begin{definition}
Let $\Cat{BD}_F(\alg{G}, \alg{T})$ be the category whose objects are triples $(Q, \sheaf{D}, f)$, where
\begin{enumerate}
\item
$Q \From \sheaf{Y} \To \ZZ$ is a Weyl-invariant quadratic form;
\item
$\sheaf{D}$ is an extension of $\sheaf{Y}$ by $\sheaf{G}_{\mult}$ (as sheaves of groups on $F_{\et}$), whose commutator is given by $\Comm(y_1, y_2)  = (-1)^{B_Q(y_1, y_2)}$;
\item
$f \From \sheaf{D}_Q \To \sheaf{D}$ is a morphism of sheaves of groups on $F_{\et}$, making the diagram \eqref{ThirdInvariantDiagram} commute.
\end{enumerate}
The morphisms in this category, from $(Q_1, \sheaf{D}_1, f_1)$ to $(Q_2, \sheaf{D}_2, f_2)$, are given as follows:  if $Q_1 \neq Q_2$ then there are no morphisms.  When $Q = Q_1 = Q_2$, the morphisms consist of morphisms $d \From \sheaf{D}_1 \To \sheaf{D}_2$ in $\Cat{CExt}_{F_{\et}}(\sheaf{Y}, \sheaf{G}_{\mult})$ which satisfy $f_2 = d \circ f_1$.
\end{definition}

For objects $(Q_1, \sheaf{D}_1, f_1)$ and $(Q_2, \sheaf{D}_2, f_2)$, define their sum to be the object $(Q, \sheaf{D}, f)$ given by:
\begin{enumerate}
\item
$Q = Q_1 + Q_2$;
\item
$\sheaf{D} = \sheaf{D}_1 \Baer \sheaf{D}_2$ (the Baer sum);
\item
Given $f_1 \From \sheaf{D}_Q \To \sheaf{D}_1$ and $f_2 \From \sheaf{D}_Q \To \sheaf{D}_2$, the universal property of the Baer sum defines a morphism $f \From \sheaf{D}_Q \To \sheaf{D} = \sheaf{D}_1 \Baer \sheaf{D}_2$.  
\end{enumerate}
This makes $\Cat{BD}_F(\alg{G}, \alg{T})$ into a Picard category.  

Since the category of extensions of $\sheaf{Y}$ by $\sheaf{G}_{\mult}$ forms a Picard stack on $F_{\et}$ (Galois descent is effective), and Galois descent is also effective for quadratic forms (global sections of the sheaf $\ssym^2(\sheaf{X})^{\sheaf{W}}$), we find a Picard stack $\stack{BD}(\alg{G}, \alg{T})$ on $F_{\et}$.  

From each object $\alg{G}'$ of the category $\Cat{CExt}_F(\alg{G}, \alg{K}_2)$, Brylinski and Deligne associate such a triple $(Q, \sheaf{D}, f)$ of invariants.  The following is a restatement of the main theorem of \cite[Theorem 7.2]{B-D}.
\begin{thm}
\label{BDMain}
The association of the three invariants, described in Section \ref{BDSection}, gives an equivalence of Picard stacks on $F_{\et}$,$$\Fun{BD} \From \stack{CExt}(\alg{G}, \alg{K}_2) \xrightarrow{\sim} \stack{BD}(\alg{G}, \alg{T}).$$
\end{thm}

When $\alg{G} = \alg{T}$ is an algebraic torus, the third invariant plays no role and the classification is simpler.  We get an equivalence of Picard stacks,
$$\Fun{BD} \From \stack{CExt}(\alg{T}, \alg{K}_2) \xrightarrow{\sim} \stack{BD}(\alg{T}).$$
The Picard category $\Cat{BD}_F(\alg{T})$ is the category of pairs $(Q, \sheaf{D})$ with $Q \From \sheaf{Y} \To \ZZ$ a quadratic form, and $\sheaf{G}_{\mult} \Into \sheaf{D} \Onto \sheaf{Y}$ a central extension with commutator $\Comm(y_1, y_2) = (-1)^{B_Q(y_1, y_2)}$.

An important special case is the following.
\begin{proposition}
The category $\Cat{CExt}_F(\alg{G}_{\mult}, \alg{K}_2)$ is equivalent to its full subcategory with objects $\{ \alg{E}_c : c \in \ZZ \}$ and morphisms the automorphisms of each $\alg{E}_c$.
\end{proposition}
\proof
If $\alg{E}$ is a central extension of $\alg{G}_{\mult}$ by $\alg{K}_2$, then by \cite[\S 3.9, 3.10]{B-D}, there exists an isomorphism in $\Cat{CExt}_F(\alg{G}_{\mult}, \alg{K}_2)$ from $\alg{E}$ to $\alg{E}_c$ for some integer $c$.

The objects $\alg{E}_c$ are not isomorphic to each other, again by \cite[\S 3.9]{B-D}, and so the only morphisms among them are the automorphisms.
\qed

\begin{example}
\label{Ec3}
The automorphisms of any object in the category $\Cat{CExt}_F(\alg{G}_{\mult}, \alg{K}_2)$ are in natural bijection with $F^\times$ by \cite[\S 3.11]{B-D}.  For the object $\alg{E}_c$, and $a \in F^\times$, the corresponding automorphism will be denoted $\alpha_a$.  Recalling that $\alg{E}_c$ is the split torsor $\alg{G}_{\mult} \times \alg{K}_2$, the automorphism $\alpha_a$ is given explicitly by
$$\alpha_a( u, \kappa) = (u, \kappa \cdot \{u,a \} ).$$
\end{example}

\subsection{The $\alg{G}_{\mult}$ analogue}

For comparison and for what comes later, we describe the much easier $\alg{G}_{\mult}$ ($=\alg{K}_1$) analogue of the main theorem \cite[Theorem 7.2]{B-D}.  For consistency in notation with later sections, we call our base field $\kk$ instead of $F$ for the moment.  Consider the category of central extensions $\Cat{CExt}_\kk(\alg{G}, \alg{G}_{\mult})$ of central extensions of $\alg{G}$ by $\alg{G}_{\mult}$ (as algebraic groups over $\kk$).  Given a maximal torus $\alg{T} \subset \alg{G}$, any such extension $\alg{G}_{\mult} \Into \alg{G}' \Onto \alg{G}$ yields an extension of tori $\alg{G}_{\mult} \Into \alg{T}' \Onto \alg{T}$ by pullback.  Taking cocharacter lattices yields
$$\ZZ \Into \sheaf{Y}' \Onto \sheaf{Y},$$
an extension of sheaves of abelian groups on $\kk_{\et}$.  

Since $\alg{G}_{\SC}$ is simply-connected semisimple, the pullback $\alg{G}_{\mult} \Into \alg{G}_{\SC}' \Onto \alg{G}_{\SC}$ splits uniquely.  Thus we find a commutative diagram of sheaves of abelian groups on $\kk_{\et}$, with exact rows.
\begin{equation}
\label{K1Diagram}
\begin{tikzcd}
\ZZ \inarrow{r} \arrow{d}{=} & \sheaf{Y}_{\SC} \times \ZZ \onarrow{r} \arrow{d}{f} & \sheaf{Y}_{\SC} \inarrow{d}{\jota} \\
\ZZ \inarrow{r} & \sheaf{Y}' \onarrow{r} & \sheaf{Y}
\end{tikzcd}
\end{equation}

Define $\Cat{EZ}_\kk(\alg{G}, \alg{T})$ to be the category of pairs $(\sheaf{Y}', f)$ where $\sheaf{Y}'$ is an object of $\Cat{Ext}(\sheaf{Y}, \ZZ)$ and $f$ is a homomorphism from the split extension $\sheaf{Y}_{\SC} \times \ZZ$ to $\sheaf{Y}'$ making the diagram \eqref{K1Diagram} commute.  A morphism from $(\sheaf{Y}_1', f_1)$ to $(\sheaf{Y}_2', f_2)$ is a morphism $\phi \From \sheaf{Y}_1' \To \sheaf{Y}_2'$ in $\Cat{Ext}(\sheaf{Y}, \ZZ)$ such that $\phi \circ f_1 = f_2$.  

The Baer sum makes both categories $\Cat{CExt}_\kk(\alg{G},\alg{G}_{\mult})$ and $\Cat{EZ}_\kk(\alg{G}, \alg{T})$ into Picard categories.  Galois descent is effective for $\alg{G}_{\mult}$-torsors just as for $\alg{K}_2$-torsors, and so we find Picard stacks $\stack{CExt}(\alg{G}, \alg{G}_{\mult})$ and $\stack{EZ}(\alg{G}, \alg{T})$ on $\kk_{\et}$.

The $\alg{G}_{\mult}$-analogue of the main result of \cite{B-D} follows.
\begin{thm}
\label{EZMain}
The construction above gives an equivalence of Picard stacks, 
$$\Fun{EZ} \From \stack{CExt}(\alg{G}, \alg{G}_{\mult}) \To \stack{EZ}(\alg{G}, \alg{T}).$$
\end{thm}
\proof
At the level of categories, this is \cite[Theorem 1.1]{Phobia}.  The compatibility with Baer sums and pullbacks is straightforward.  
\qed

\section{Integral models}
\label{IMSection}
In this section, let $\OO$ be a Dedekind domain with fraction field $F$.  Write $S = \Spec(\OO)$ and $S^{(1)}$ for the closed points of $S$.  Write $\eta \From \Spec(F) \Into S$ for the inclusion of the generic point.  When $s \in S^{(1)}$, write $\OO_s$ for the local ring at $s$, and $\kk(s)$ for its residue field.

Let $\alg{G}$ be a reductive group over $\OO$, with maximal torus $\alg{T}$ over $\OO$, and generic fibre $\alg{G}_F$.  For all $s \in S^{(1)}$, write $\alg{\bar G}_s$ for the resulting reductive group over the residue field $\kk(s)$.

\begin{remark}
When $\OO$ is the ring of integers in a nonarchimedean local field, Tits (\cite[\S 3.8.1]{Tits}) demonstrates that $\alg{G}_F$ is a quasisplit reductive group over $F$, $\alg{G}_F$ splits over an unramified extension of $F$, and $\alg{G}(\OO)$ is a hyperspecial maximal compact subgroup of $\alg{G}(F)$.
\end{remark}

Write $\Cat{CExt}_\OO(\alg{G}, \alg{K}_2)$ for the Picard category of central extensions of sheaves of groups on $\OO_{\Zar}$, as studied in \cite{B-D}.  Pulling back via $\eta \From \Spec(F) \Into S$ yields an additive functor of Picard categories,
$$\eta^\ast \From \Cat{CExt}_\OO(\alg{G}, \alg{K}_2) \To \Cat{CExt}_F(\alg{G}_F, \alg{K}_2).$$

\begin{definition}
When $\alg{G}_F'$ is a central extension of $\alg{G}_F$ by $\alg{K}_2$, an \defined{integral model} (or an \defined{$\OO$-model}) of $\alg{G}_F'$  is a pair $(\alg{G}', \iota)$ with $\alg{G}'$ a central extension of $\alg{G}$ by $\alg{K}_2$, and $\iota \From \eta^\ast \alg{G}' \To \alg{G}_F'$ an isomorphism in $\Cat{CExt}_F(\alg{G}_F, \alg{K}_2)$.  
\end{definition}

The key result for understanding central extensions of $\alg{G}$ by $\alg{K}_2$ is the Quillen-Gersten resolution, discussed in a broader setting here.

\subsection{Quillen-Gersten resolution}

When $\alg{X}$ is a smooth scheme of finite type over $S$, let $\alg{X}^{(i)}$ be the set of points of codimension $i$ in $\alg{X}$.  The Quillen-Gersten complex is composed of terms
$$\alg{Q}_n^i = \bigoplus_{x \in \alg{X}^{(i)}} \iota_{x \ast} \alg{K}_{n-i}(\kk(x)), \quad (0 \leq i \leq n)$$
with $\kk(x)$ the appropriate residue fields, and morphisms 
$$0 \To \alg{K}_n \To \alg{Q}_n^0 \To \alg{Q}_n^1 \To \cdots \To \alg{Q}_n^n \To 0$$
given by residue maps in $\alg{K}$-theory.  These $\alg{Q}_n^i$ are flabby sheaves on $\alg{X}_{\zar}$.  When $n = 0$, $\alg{K}_0 = \alg{Q}_0^0 = \ZZ$, the constant sheaf.  When $n = 1$, we get a familiar resolution of the Zariski sheaf $\alg{K}_1 = \alg{G}_{\mult}$, whose cohomology relates Weil divisors to Cartier divisors,
$$0 \To \alg{K}_1 \To \alg{Q}_1^0 \To \alg{Q}_1^1 \To 0.$$

While the exactness of the Quillen-Gersten complex remains an open problem (called \defined{Gersten's conjecture}), the following results for $\alg{K}_2$ are known.
\begin{thm}
\label{GerstenCases}
Assume one of the following:
\begin{enumerate}
\item
$\alg{X}$ is a smooth scheme of finite type over a field;
\item
$\alg{X}$ is a smooth scheme of finite type over $S = \Spec(\OO)$, with $\OO$ a discrete valuation ring;
\item
$\alg{X} = S = \Spec(\OO)$ or $\alg{X}$ an open $\OO$-subscheme of $\Spec(\OO[\form])$, with $\OO$ a Dedekind domain with finite residue fields.
\end{enumerate}
Then the complex of sheaves on $\alg{X}_{\zar}$
\begin{equation}
\label{GQComplex}
0 \To \alg{K}_2 \To \alg{Q}_2^0 \To \alg{Q}_2^1 \To \alg{Q}_2^2 \To 0
\end{equation}
is exact.
\end{thm}
\proof
The first case is a (special case) of a result of Quillen \cite{Quillen}.  The second is a result of Bloch \cite[Corollary A.2]{Blo}. The third case is a direct consequence of \cite{ShermanGlobal}.  
\qed

We say that Gersten's conjecture holds in \defined{weight two} for $\alg{X} \To S$ if the complex \eqref{GQComplex} is an exact sequence of sheaves on $\alg{X}_{\zar}$.  The following two results are applications of the Quillen-Gersten resolution, based on work of Sherman \cite{ShermanGlobal} and the Quillen localization sequence.
\begin{proposition}
\label{KCohomology}
Suppose that $\OO$ is a Dedekind domain with fraction field $F$, $S = \Spec(\OO)$, and at least one of the following holds:
\begin{enumerate}
\item
$\OO$ contains a field, or...
\item
All residue fields $\kk(s)$ of $\OO$ are finite.
\end{enumerate}
Then we have
$$H_{\zar}^0(S, \alg{K}_2) = \Im \left( K_2(\OO) \To K_2(F) \right), \quad H_{\zar}^1(S, \alg{K}_2) = \Ker(K_1(\OO) \To K_1(F));$$ 
$$H_{\zar}^0(S, \alg{K}_1) = \Im \left( K_1(\OO) \To K_1(F) \right), \quad H_{\zar}^1(S, \alg{K}_1) = \Ker(K_0(\OO) \To K_0(F));$$
$$H_{\zar}^0(\alg{G}_{\mult/S}, \alg{K}_2) = H_{\zar}^0(S, \alg{K}_2) \oplus H_{\zar}^0(S, \alg{K}_1);$$
$$H_{\zar}^1(\alg{G}_{\mult/S}, \alg{K}_2) = H_{\zar}^1(S, \alg{K}_2) \oplus H_{\zar}^1(S, \alg{K}_1).$$
\end{proposition}
\proof
This follows directly from \cite[Corollaries 2.5,2.6]{ShermanGlobal} and \cite[Theorem 4.3]{She}
\qed

Under the hypotheses of this corollary, we find a residue homomorphism
$$\partial \From H_{\zar}^0(\alg{G}_{\mult/S}, \alg{K}_2) \To H_{\zar}^0(S, \alg{K}_1).$$
In cases of arithmetic interest, we find the vanishing of some cohomology groups.
\begin{corollary}
\label{H1KVan}
Suppose that $S$ is one of the following base schemes:
\begin{enumerate}
\item
$S = \Spec(\OO)$ with $\OO$ a discrete valuation ring;
\item
$S = \Spec(\OO)$, with $\OO = \OO_\SS$ the ring of $\SS$-integers in a global field, and $\SS$ sufficiently large so that $\Pic(\OO) = 0$.
\end{enumerate}
Then
$$H_{\zar}^1(\alg{G}_{\mult/S}, \alg{K}_2) = H_{\zar}^1 (S, \alg{K}_2) = 0.$$
\end{corollary}
\proof
We apply the formulae of Proposition \ref{KCohomology} throughout.  For both types of base scheme, we have 
$$\Pic(S) = \Ker \left( K_0(\OO) \To K_0(F) \right) = 0.$$
Hence in both cases we have $H_{\zar}^1(S, \alg{K}_1) = 0$.  It follows that
$$H_{\zar}^1(\alg{G}_{\mult/S}, \alg{K}_2) = H_{\zar}^1(S, \alg{K}_2) = \Ker \left( K_1(\OO) \To K_1(F) \right).$$
When $\OO$ is a discrete valuation ring, $K_1(\OO) = \OO^\times$ and the map $K_1(\OO) \To K_1(F)$ is just the inclusion $\OO^\times \Into F^\times$.  When $\OO = \OO_\SS$, the ring of $\SS$-integers in a global field, $K_1(\OO) = \OO^\times$ by the main theorem \cite[Theorem 4.1]{BMS} of Bass, Milnor, and Serre.

Hence in both cases, we have
$$H_{\zar}^1(\alg{G}_{\mult/S}, \alg{K}_2) = H_{\zar}^1(S, \alg{K}_2) = 0.$$
\qed

\subsection{Integral models and residual extensions}

Return to the general setting where $\OO$ is a Dedekind domain with fraction field $F$, and $S = \Spec(\OO)$.  A reductive group $\alg{G}$ over $S$ will be fixed.  Write $j \From \alg{G}_F \Into \alg{G}$ for the inclusion of the general fibre.  For each $s \in S^{(1)}$, let $i_s \From \alg{\bar G}_s \Into \alg{G}$ be the inclusion of the corresponding special fiber.  \textbf{Assuming Gersten's conjecture} (in weight two for finite-type schemes over $S$), and following the local results mentioned in \cite[Remarks 12.14 (iii)]{B-D}, the Quillen-Gersten resolution of $\alg{K}_2$ gives a short exact sequence of sheaves on $\alg{G}_{\zar}$,
\begin{equation}
\label{K2Seq}
\alg{K}_2 \Into j_\ast \alg{K}_2 \Onto \bigoplus_{s \in S^{(1)}} i_{s \ast} \alg{K}_1.
\end{equation}
To see this, we write down the Quillen-Gersten resolution of each term.  We partition the points of the scheme $\alg{G}$ according to whether they lie over the generic point $\Spec(F)$ of $S$ or over a closed point $s \in S^{(1)}$.  Write $g_j$ for the generic point of $\alg{G}_F$ and $g_s$ for the generic point of the special fibre $\alg{\bar G}_s$.  This gives
$$\alg{G}^{(0)} = \{ g_j \}, \quad \alg{G}^{(1)} = \alg{G}_F^{(1)} \sqcup \{ g_s : s \in S^{(1)} \},$$
$$\alg{G}^{(2)} = \alg{G}_F^{(2)} \sqcup \bigsqcup_{s \in S^{(1)}} \alg{\bar G}_s^{(1)}.$$

We decompose the Quillen-Gersten resolution according to these partitions, and abbreviate by writing $\iota_\ast$ for every pushforward of a sheaf from a point to the scheme $\alg{G}$.
$$\begin{tikzcd}[column sep = 1.2em]
\alg{K}_2 \inarrow{d} & j_\ast \alg{K}_2 \inarrow{d} & \phantom{a} \\ 
\iota_\ast \alg{K}_2(\kk(g_j)) \arrow{d} &   \iota_\ast \alg{K}_2(\kk(g_j)) \arrow{d} &  \bigoplus\limits_{s \in S^{(1)}} \iota_{s \ast} \alg{K}_1 \inarrow{d}  \\
\bigoplus\limits_{s \in S^{(1)}} \iota_\ast \alg{K}_1(\kk(g_s)) \oplus \bigoplus\limits_{x \in \alg{G}_F^{(1)}} \iota_\ast  \alg{K}_1(\kk(x)) \arrow{d} &  \bigoplus\limits_{x \in \alg{G}_F^{(1)}} \iota_\ast  \alg{K}_1(\kk(x)) \arrow{d} & \bigoplus\limits_{s \in S^{(1)}} \iota_\ast \alg{K}_1(\kk(g_s)) \arrow{d} \\
\bigoplus\limits_{\substack{s \in S^{(1)} \\ y \in \alg{\bar G}_s^{(1)}}} \iota_\ast \ZZ \oplus \bigoplus\limits_{x \in \alg{G}_F^{(2)}} \iota_\ast \ZZ & \bigoplus\limits_{x \in \alg{G}_F^{(2)}} \iota_\ast \ZZ  & \bigoplus\limits_{\substack{s \in S^{(1)} \\ y \in \alg{\bar G}_s^{(1)}}} \iota_\ast \ZZ 
\end{tikzcd}$$
A diagram chase yields the desired short exact sequence \eqref{K2Seq} of sheaves of abelian groups on $\alg{G}_{\zar}$.

The short exact sequence \eqref{K2Seq} of sheaves gives an exact sequence of Picard categories and additive functors by Proposition \ref{LExact},
\begin{equation}
\label{CExtExact}
0 \To \Cat{CExt}_\OO(\alg{G}, \alg{K}_2) \To \Cat{CExt}_\OO(\alg{G}, j_\ast \alg{K}_2) \To \Cat{CExt}_\OO \left( \alg{G},  \bigoplus_{s \in S^{(1)}} i_{s \ast} \alg{K}_1 \right).
\end{equation}

Adjunction identifies the Picard category $\Cat{CExt}_\OO(\alg{G}, j_\ast \alg{K}_2)$ with the Picard category $\Cat{CExt}_F(\alg{G}_F, \alg{K}_2)$.  Similarly, for every $s \in S^{(1)}$, we identify the Picard categories $\Cat{CExt}_\OO(\alg{G}, i_{s \ast} \alg{K}_1)$ with $\Cat{CExt}_{\kk(s)}(\alg{\bar G}_s, \alg{\bar G}_{\mult})$.  We find an exact sequence of Picard categories and additive functors,
\begin{equation}
\label{CExtExactBetter}
0 \To \Cat{CExt}_\OO(\alg{G}, \alg{K}_2) \xrightarrow{\eta^\ast} \Cat{CExt}_F(\alg{G}_F, \alg{K}_2) \xrightarrow{\partial} \bigoplus_{s \in S^{(1)}} \Cat{CExt}_{\kk(s)}(\alg{\bar G}_s, \alg{\bar G}_{\mult}).
\end{equation}
The functor $\eta^\ast$ is the pullback via $\eta \From \Spec(F) \Into S$.  The functor $\partial$ is described in more detail in the construction of \cite[\S 12.11]{B-D}; it can be described as a direct sum of functors $\partial = \bigoplus_{s \in S^{(1)}} \partial_s$.  When $\alg{G}_F' \in \Cat{CExt}(\alg{G}_F, \alg{K}_2)$, the functor $\partial_s$ provides an extension $\alg{\bar G}_s' = \partial_s \alg{G}_F' \in \Cat{CExt}_{\kk(s)}(\alg{\bar G}_s, \alg{\bar G}_{\mult})$.  We call this $\alg{\bar G}_s'$ the \defined{residual extension} of $\alg{G}_F'$ at $s$, and it is a focus of the next section.

Exactness for the sequence of Picard categories gives the following result.
\begin{proposition}
(Assume Gersten's conjecture in weight two for finite-type schemes over $\OO$.)
\label{ExistIntModel}
The functors $\eta^\ast$ and $\partial$ provide an equivalence of Picard categories from
\begin{itemize}
\item
the category $\Cat{CExt}_\OO(\alg{G}, \alg{K}_2)$ of central extensions defined over $\OO$, to...
\item
the category of pairs $(\alg{G}_F', \tau)$ where $\alg{G}_F' \in \Cat{CExt}_F(\alg{G}_F, \alg{K}_2)$, and $\tau = (\tau_s)_{s \in S^{(1)}}$ is a family of splittings of every residual extension $\alg{\bar G}_{\mult} \Into \alg{\bar G}_s' \Onto \alg{\bar G}_s$.
\end{itemize}
\end{proposition}

\begin{corollary}
\label{UniqueIntModel}
(Assume Gersten's conjecture in weight two for finite-type schemes over $\OO$.)  Suppose that $\alg{G}$ is a reductive group over $S$, with simply-connected semisimple geometric fibres.  Let $\alg{G}_F'$ be an extension of $\alg{G}_F$ by $\alg{K}_2$ defined over $F$.  Then there is a unique, up to unique isomorphism, $\OO$-model $(\alg{G}', \iota)$ of $\alg{G}_F$.
\end{corollary}
\proof
In this setting, every residual extension
$$\alg{\bar G}_{\mult} \Into \alg{\bar G}_s' \Onto \alg{\bar G}_s$$
is an extension of a simply-connected semisimple group $\alg{\bar G}_s$ by $\alg{\bar G}_{\mult}$, defined over $\kk(s)$.  Such an extension splits uniquely, yielding the desired result.
\qed

\begin{example}
\label{Ec4}
Suppose that $\OO$ is a discrete valuation ring with residue field $\kk$ and fraction field $F$.  In particular, Gersten's conjecture holds in weight two for finite-type schemes over $\OO$, by \cite{Blo}.  Applying the functor $\partial$ to the object $\alg{E}_c \in \Cat{CExt}_F(\alg{G}_{\mult/F}, \alg{K}_2)$ (of Example \ref{Ec1}) yields an extension $\alg{\bar E}_c \defeq \partial \alg{E}_c \in \Cat{Ext}_\kk(\alg{\bar G}_{\mult}, \alg{\bar G}_{\mult})$.  

The triviality of the torsor $\alg{E}_c = \alg{G}_{\mult} \times \alg{K}_2$ trivializes the $\alg{\bar G}_{\mult}$-torsor, so we have $\alg{\bar E}_c = \alg{\bar G}_{\mult} \times \alg{\bar G}_{\mult}$ (as a $\alg{\bar G}_{\mult}$-torsor on $\alg{\bar G}_{\mult}$).  Multiplication in $\alg{E}_c$, given by $(u, \alpha) \cdot (v, \beta) = (uv, \alpha \beta \cdot \{ u,v \}^c)$ yields multiplication in $\alg{E}_c$:
$$(\bar u, \bar a) \cdot (\bar v, \bar b) = (\bar u \bar v, \bar a \bar b \cdot \partial \{ u,v \} ),$$
for all $u,v \in \OO^\times$ projecting to $\bar u, \bar v \in \kk^\times$.  But the residue symbol $K_2(F) \xrightarrow{\partial} K_1(\kk)$ is trivial on the image of $K_2(\OO)$, and so we have
$$(\bar u, \bar a) \cdot (\bar v, \bar b) = (\bar u \bar v, \bar a \bar b).$$
In other words, the residual extension $\alg{\bar E}_c$ equals the split extension $\alg{\bar G}_{\mult} \times \alg{\bar G}_{\mult}$.  This endows $\alg{E}_c$ with an $\OO$-model. 

One may similarly trace through the automorphism $\alpha_a \in \Aut(\alg{E}_c)$ (see Example \ref{Ec3}) to find that
$$\partial \alpha_a (\bar u, \bar v) = (\bar u, \bar v \cdot \bar u^{\val a})$$
for all $(\bar u, \bar v) \in \alg{\bar E}_c = \alg{\bar G}_{\mult} \times \alg{\bar G}_{\mult}$.  
\end{example}

\subsection{Three invariants for integral models}
\label{InvariantsIntegral}

We may formulate the three invariants of Brylinski-Deligne, working with reductive groups $\alg{G}$ over a Dedekind domain $\OO$ of one of the following sorts:
\begin{enumerate}
\item
$\OO$ a DVR with finite residue field;
\item
$\OO$ a DVR containing a field;
\item
 $\OO = \OO_\SS$ the ring of $\SS$-integers in a global field of characteristic $p$, with $\SS$ sufficiently large so that $\alg{Pic}(\OO) = 0$;
\item
$\OO = \OO_\SS$ the ring of $\SS$-integers in a function field, with $\SS$ sufficiently large so that $\alg{Pic}(\OO) = 0$, \textbf{assuming Gersten's conjecture} holds in weight two for smooth schemes of finite type over $\OO$.
\end{enumerate}

Fix a reductive group $\alg{G}$ over such a Dedekind domain $\OO$, with a maximal torus $\alg{T}$ over $\OO$, as well as an extension $\alg{K}_2 \Into \alg{G}' \Onto \alg{G}$ in $\Cat{CExt}_\OO(\alg{G}, \alg{K}_2)$.  We write $\eta \From \Spec(F) \Into S = \Spec(\OO)$ in what follows.  
\subsubsection{First invariant}
The sheaf $\sheaf{Y}$ of cocharacters of $\alg{T}$ is a local system of abelian groups on $\OO_{\et}$.  Similarly, the sheaf $\sheaf{Y}_F$ of cocharacters of $\alg{T}_F$ is a local system of abelian groups on $F_{\et}$.  

Since $\OO$ is a Dedekind domain, the pullback functor which sends a local system of abelian groups on $\OO_{\et}$ to the resulting local system on $F_{\et}$ is {\em fully faithful}.  Indeed, this corresponds to the fact that the map of \'etale fundamental groups $\pi_1^{\et}(\Spec(F)) \To \pi_1^{\et}(\Spec(\OO))$ is surjective, with respect to any geometric base point $\bar F / F$.  

It follows that the first invariant of $\alg{G}_F' \in \Cat{CExt}_F(\alg{G}_F, \alg{K}_2)$, a global section
$$Q_F \in H^0 \left( F_{\et}, \ssym^2(\sheaf{X}_F)^{\sheaf{W}_F} \right)$$
arises from a unique global section,
$$Q \in H^0 \left( \OO_{\et}, \ssym^2(\sheaf{X})^{\sheaf{W}} \right)$$

In this way, the first invariant $Q_F$ for the object $\alg{G}_F' \in \Cat{CExt}_F(\alg{G}_F, \alg{K}_2)$ yields a \defined{first invariant} $Q$ for the object $\alg{G}' \in \Cat{CExt}_\OO(\alg{G}, \alg{K}_2)$.  

\begin{remark}
This construction of the first invariant relies only on the assumption that $S$ is the spectrum of a Dedekind domain.
\end{remark}

\subsubsection{Second invariant}

The construction of the second invariant, the extension of sheaves $\sheaf{G}_{\mult} \Into \sheaf{D} \Onto \sheaf{Y}$, works as well over $\OO_{\et}$ as it does over $F_{\et}$.  Step by step, we begin with $\alg{K}_2 \Into \alg{T}' \Onto \alg{T}$.  We work locally on $\OO_{\et}$, beginning with a connected finite \'etale $U \To S = \Spec(\OO)$ for which $\sheaf{Y}_U$ is constant.  Take global sections over $\alg{G}_{\mult/U}$ to obtain a short exact sequence
\begin{equation}
\label{SESOverU}
H_{\zar}^0(\alg{G}_{\mult/U}, \alg{K}_2) \Into H_{\zar}^0(\alg{G}_{\mult/U}, \alg{T}') \To H_{\zar}^0(\alg{G}_{\mult/U}, \alg{T})
\end{equation}
Here we apply Corollary \ref{H1KVan} for vanishing of $H_{\zar}^1(\alg{G}_{\mult/U}, \alg{K}_2)$.   Proposition \ref{KCohomology} provides a homomorphism 
$$\partial \From H_{\zar}^0(\alg{G}_{\mult/U}, \alg{K}_2) \To H_{\zar}^0(U, \alg{K}_1).$$
In the settings under consideration, we have
$$H_{\zar}^0(U, \alg{K}_1) = \Im \left( K_1(U) \To K_1(F) \right) = \sheaf{G}_{\mult}[U].$$  
Thus, pushing out \eqref{SESOverU} via $\partial$ gives a short exact sequence
$$\sheaf{G}_{\mult}[U] \Into \partial_\ast H_{\zar}^0(\alg{G}_{\mult/U}, \alg{T}') \Onto H_{\zar}^0(\alg{G}_{\mult/U}, \alg{T}).$$
Pulling back via the canonical homomorphism $h \From \sheaf{Y}[U] \To H_{\zar}^0(\alg{G}_{\mult/U}, \alg{T})$ defines $\sheaf{D}^\circ[U] = h^\ast \partial_\ast H_{\zar}^0(\alg{G}_{\mult/U}, \alg{T}')$, lying in a short exact sequence
$$\sheaf{G}_{\mult}[U] \Into \sheaf{D}^\circ[U] \Onto \sheaf{Y}[U].$$
The construction is functorial in $U \To S$, and thus defines a short exact sequence of sheaves of groups on $\OO_{\et}$,
$$\sheaf{G}_{\mult} \Into \sheaf{D}^\circ \Onto \sheaf{Y}.$$
This defines an additive functor of Picard categories,
$$\Cat{CExt}_\OO(\alg{G}, \alg{K}_2) \To \Cat{CExt}_{\OO_{\et}}(\sheaf{Y}, \sheaf{G}_{\mult}).$$
\begin{remark}
A priori, the functoriality of the construction in $U \To S$ gives a short exact sequence of {\em presheaves} on $\OO_{\et}$.  But since $\sheaf{Y}$ and $\sheaf{G}_{\mult}$ satisfy the sheaf axiom, the presheaf $\sheaf{D}^\circ$ also satisfies the sheaf axiom.
\end{remark}

\begin{remark}
This construction of the second invariant relies on all of our assumptions on $\OO$, except for the requirement of Gersten's conjecture.
\end{remark}

\subsubsection{Third invariant}
The construction of the third invariant relies on the first two invariants, and the uniqueness, up to unique isomorphism, of an extension $\alg{G}_{\SC}'$ with first invariant equal to a Weyl-invariant quadratic form $Q$.  The uniqueness, up to unique isomorphism, of an extension $\alg{K}_2 \Into \alg{G}_{\SC, F}' \Onto \alg{G}_{\SC,F}$ with first invariant $Q$ is guaranteed by Theorem \ref{UniqueCoversSCSS}.  The uniqueness, up to unique isomorphism, of an $\OO$-model of $\alg{G}_{\SC,F}'$ is guaranteed by Corollary \ref{UniqueIntModel} (whose proof relied on Gersten's conjecture).  Define $\sheaf{D}_Q^\circ$ to be the second invariant of the unique $\OO$-model of the canonical extension $\alg{G}_{\SC,Q}'$.

We find a commutative diagram of sheaves of groups on $\OO_{\et}$.
\begin{equation}
\label{IntDDiagram}
\begin{tikzcd}
\sheaf{G}_{\mult} \inarrow{r} \arrow{d}{=} & \sheaf{D}_Q^\circ \onarrow{r} \inarrow{d}{f^\circ} & \sheaf{Y}_{\SC} \inarrow{d}{\jota} \\
\sheaf{G}_{\mult} \inarrow{r} & \sheaf{D}^\circ \onarrow{r} & \sheaf{Y}
\end{tikzcd}
\end{equation}
This gives the third invariant, a homomorphism $f^\circ$ of sheaves on $\OO_{\et}$.

\subsubsection{Category of triples}

Define $\Cat{BD}_\OO(\alg{G}, \alg{T})$ be the category whose objects are triples $(Q, \sheaf{D}^\circ, f^\circ)$, where
\begin{enumerate}
\item
$Q \in H^0 \left( \OO_{\et}, \ssym^2(\sheaf{X})^{\sheaf{W}} \right)$ is a Weyl-invariant quadratic form on $\sheaf{Y}$;
\item
$\sheaf{D}^\circ \in \Cat{CExt}_{\OO_{\et}}(\sheaf{Y},\sheaf{G}_{\mult})$ whose commutator is given by
$$\Comm(y_1, y_2) = (-1)^{B_Q(y_1, y_2)}.$$
\item
$f^\circ \From \sheaf{D}_Q^\circ \To \sheaf{D}^\circ$ is a morphism of sheaves of groups on $\OO_{\et}$ making the diagram \eqref{IntDDiagram} commute.
\end{enumerate}

The three invariants give an additive functor of Picard categories,
$$\Fun{BD}_\OO \From \Cat{CExt}_\OO(\alg{G}, \alg{K}_2) \To \Cat{BD}_\OO(\alg{G}, \alg{T}).$$
Pulling back via $\eta \From \Spec(F) \To \Spec(\OO)$ yields a diagram of Picard categories and additive functors, which commutes up to natural isomorphism (i.e., the diagram 2-commutes).
\begin{equation}
\label{BDPull}
\begin{tikzcd}
\Cat{CExt}_\OO(\alg{G}, \alg{K}_2) \arrow{r}{\eta^\ast} \arrow{d}{\Fun{BD}_\OO} & \Cat{CExt}_F(\alg{G}_F, \alg{K}_2) \arrow{d}{\Fun{BD}_F}\\
\Cat{BD}_\OO(\alg{G}, \alg{T}) \arrow{r}{\eta^\ast} & \Cat{BD}_F(\alg{G}_F, \alg{T}_F) 
\end{tikzcd}
\end{equation}

\section{The residual extension}
\label{RESection}

Throughout this section, we focus on the local case, where $\OO$ is a DVR with residue field $\kk$ and fraction field $F$.  We assume that $\OO$ contains a field, or that the residue field $\kk$ is finite, so that the results of the previous section hold.  Write $\eta \From \Spec(F) \Into \Spec(\OO)$, and $\iota \From \Spec(\kk) \Into \Spec(\OO)$ for the inclusions of generic point and closed point into the scheme $\Spec(\OO)$.

Suppose that $\alg{G}$ is a reductive group over $\OO$ with maximal torus $\alg{T}$ over $\OO$.  Their special fibres are denoted $\alg{\bar G}$ and $\alg{\bar T}$; the latter is a maximal $\kk$-torus in the former.  In this setting, the residual extension is given by a functor
$$\partial \From \Cat{CExt}_F(\alg{G}_F, \alg{K}_2) \To \Cat{CExt}_{\kk}(\alg{\bar G}, \alg{\bar G}_{\mult}).$$

Theorems \ref{BDMain} and \ref{EZMain} give equivalences of Picard categories 
$$\Fun{BD}_F \From \Cat{CExt}_F(\alg{G}_F, \alg{K}_2) \xrightarrow{\sim} \Cat{BD}_F(\alg{G}_F, \alg{T}_F), \quad \Fun{EZ}_{\kk} \From \Cat{CExt}_{\kk}(\alg{\bar G}, \alg{\bar G}_{\mult}) \xrightarrow{\sim} \Cat{EZ}_\kk(\alg{\bar G}, \alg{\bar T}).$$
We begin by describing an additive functor $\Fun{val} \From \Cat{BD}_F(\alg{G}_F, \alg{T}_F) \To \Cat{EZ}_\kk(\alg{\bar G}, \alg{\bar T})$, and a natural isomorphism $N \From \Fun{val} \circ \Fun{BD}_F \xRightarrow{\sim} \Fun{EZ}_{\kk} \circ \partial$ which makes the following diagram of Picard categories and additive functors 2-commute.
$$\begin{tikzcd}
\Cat{CExt}_F(\alg{G}_F, \alg{K}_2) \arrow{r}{\partial} \arrow{d}{\Fun{BD}_F} & \Cat{CExt}_{\kk}(\alg{\bar G}, \alg{\bar G}_{\mult}) \arrow{d}{\Fun{EZ}_{\kk}} \\
\Cat{BD}_F(\alg{G}_F, \alg{T}_F) \arrow{r}{\Fun{val}} \arrow[Rightarrow, bend right=10]{ur}{N} & \Cat{EZ}_\kk(\alg{\bar G}, \alg{\bar T})
\end{tikzcd}$$

This answers Question 12.13(i) of \cite{B-D}, in the hyperspecial case.

\subsection{The valuation functor}
\subsubsection{Tori}
As $\alg{T}$ is a torus over $\OO$, the cocharacter lattice of $\alg{T}$ is a local system $\sheaf{Y}$ on $\OO_{\et}$.  Here we define the \defined{valuation functor},
$$\Fun{val} \From \Cat{BD}_F(\alg{T}_F) \To \Cat{EZ}_\kk(\alg{\bar T}).$$
An object of $\Cat{BD}_F(\alg{T}_F)$ is a pair $(Q, \sheaf{D})$ with $Q \From \eta^\ast \sheaf{Y} \To \ZZ$ a quadratic form, and $\sheaf{D} \in \Cat{CExt}_{F_{\et}}(\eta^\ast \sheaf{Y}, \sheaf{G}_{\mult/F})$.  The category $\Cat{EZ}_\kk(\alg{\bar T})$ is the category $\Cat{Ext}_{\kk_{\et}}(\iota^\ast \sheaf{Y}, \ZZ)$.  Note that pulling back $\sheaf{Y} \mapsto \eta^\ast \sheaf{Y}$ or $\sheaf{Y} \mapsto \iota^\ast \sheaf{Y}$ does not lose information.  

The functor $\Fun{val}$ proceeds through forgetting $Q$, and applying a functor
$$\Fun{val} \From \Cat{CExt}_{F_{\et}}(\eta^\ast \sheaf{Y}, \sheaf{G}_{\mult/F}) \To \Cat{Ext}_{\kk_{\et}}(\iota^\ast \sheaf{Y}, \ZZ).$$

In the case of a split torus, we can do away with sheaves, and writing $D = \sheaf{D}[F]$, the functor $\Fun{val}$ is the pushout $\val_\ast$ described by the diagram below.
\begin{equation}
\label{valsplit}
\begin{tikzcd}
F^\times \inarrow{r} \arrow{d}{\val} & D \arrow{d} \onarrow{r} & Y \arrow{d}{=}  \\
\ZZ \inarrow{r} & \val_\ast(D) \onarrow{r} & Y
\end{tikzcd}
\end{equation}
Generally, we exploit the short exact sequence of sheaves of abelian groups on $\OO_{\et}$,
$$\sheaf{G}_{\mult} \Into \eta_\ast \sheaf{G}_{\mult/F} \xtwoheadrightarrow{\val} \iota_\ast \ZZ.$$
This defines, by Proposition \ref{LExact}, a left-exact sequence of Picard categories and additive functors,
$$0 \To \Cat{CExt}_{\OO_{\et}}(\sheaf{Y}, \sheaf{G}_{\mult}) \To \Cat{CExt}_{\OO_{\et}}(\sheaf{Y}, \eta_\ast \sheaf{G}_{\mult/F}) \xrightarrow{\Fun{val}} \Cat{CExt}_{\OO_{\et}}(\sheaf{Y}, \iota_\ast \ZZ).$$

Pullbacks and pushouts provide equivalences of Picard categories,
$$\Cat{CExt}_{\OO_{\et}}(\sheaf{Y}, \eta_\ast \sheaf{G}_{\mult/F}) \xleftrightarrow{\sim} \Cat{CExt}_{F_{\et}}(\eta^\ast \sheaf{Y}, \sheaf{G}_{\mult/F}).$$
$$\Cat{CExt}_{\OO_{\et}}(\sheaf{Y}, \iota_\ast \ZZ) \xleftrightarrow{\sim} \Cat{CExt}_{\kk_{\et}}(\iota^\ast \sheaf{Y}, \ZZ)$$

Assembling these gives an exact sequence of Picard categories and exact functors,
$$0 \To \Cat{CExt}_{\OO_{\et}}(\sheaf{Y}, \sheaf{G}_{\mult}) \To \Cat{CExt}_{F_{\et}}(\eta^\ast \sheaf{Y}, \sheaf{G}_{\mult/F}) \xrightarrow{\Fun{val}} \Cat{CExt}_{\kk_{\et}}(\iota^\ast \sheaf{Y}, \ZZ).$$
Adding a quadratic form to the mix gives an important result.
\begin{thm}
\label{ExactBD}
The following sequence of Picard categories and additive functors is exact.
$$\begin{tikzcd}[row sep = 1em]
0 \arrow{r} & \Cat{BD}_\OO(\alg{T}) \arrow{r}{\eta^\ast} & \Cat{BD}_F(\alg{T}_F) \arrow{r}{\Fun{val}} & \Cat{EZ}_\kk(\alg{\bar T}) \\
\phantom{0} & (Q, \sheaf{D}^\circ) \arrow[mapsto]{r} & (Q, \eta^\ast \sheaf{D}^\circ) & \phantom{0} \\
\phantom{0} & \phantom{0} & (Q, \sheaf{D}) \arrow[mapsto]{r} & \Fun{val} ( \sheaf{D} ) 
\end{tikzcd}$$
\end{thm}

The functors $\eta^\ast$ and $\Fun{val}$ are compatible with pullbacks; if $\tau \From \alg{T}_0 \To \alg{T}$ is a morphism of tori over $\OO$, corresponding to a morphism of cocharacter lattices $\tau \From \sheaf{Y}_0 \To \sheaf{Y}$, then pulling back defines functors $\tau^\ast$ fitting into a diagram of Picard categories and additive functors.
$$\begin{tikzcd}
0 \arrow{r} & \Cat{BD}_\OO(\alg{T}) \arrow{r}{\eta^\ast} \arrow{d}{\tau^\ast} & \Cat{BD}_F(\alg{T}_F) \arrow{r}{\Fun{val}} \arrow{d}{\tau^\ast} & \Cat{EZ}_\kk(\alg{\bar T}) \arrow{d}{\tau^\ast} \\
0 \arrow{r} & \Cat{BD}_\OO(\alg{T}_0) \arrow{r}{\eta^\ast} & \Cat{BD}_F(\alg{T}_{0,F}) \arrow{r}{\Fun{val}} & \Cat{EZ}_\kk(\alg{\bar T}_0)
\end{tikzcd}$$
The natural isomorphisms which express the commutativity of pulling back with pushing out make this diagram 2-commute.
\subsubsection{Reductive groups}
Now we consider the reductive group $\alg{G}$ over $\OO$ with maximal torus $\alg{T}$ over $\OO$.  An object of $\Cat{BD}_F(\alg{G}_F, \alg{T}_F)$ is a triple $(Q, \sheaf{D}, f \From \sheaf{D}_Q \To \sheaf{D})$.  Applying the functor $\Fun{val}$ and its compatibility with pullback yields a commutative diagram with exact rows,
$$\begin{tikzcd}
\ZZ \inarrow{r} \arrow{d}{=} & \Fun{val}(\sheaf{D}_Q) \onarrow{r} \arrow{d}{\Fun{val} f} & \sheaf{Y}_{\SC} \arrow{d}{\jota} \\
\ZZ \inarrow{r} & \Fun{val}(\sheaf{D}) \onarrow{r} & \sheaf{Y}
\end{tikzcd}$$
On the other hand, the 2-commutativity of the diagram
$$\begin{tikzcd}
\alg{G}_Q' \in \Cat{CExt}_\OO(\alg{G}_{\SC}, \alg{K}_2) \arrow{r}{\eta^\ast} \arrow{d}{\Fun{BD}} & \Cat{CExt}_F(\alg{G}_{\SC,F}, \alg{K}_2) \ni \alg{G}_{Q,F}' \arrow{d}{\Fun{BD}} \\
(Q, \sheaf{D}_Q^\circ) \in \Cat{BD}_\OO(\alg{G}_{\SC}, \alg{T}_{\SC}) \arrow{r}{\eta^\ast} & \Cat{BD}_F(\alg{T}_{\SC,F}, \alg{K}_2) \ni (Q, \sheaf{D}_Q)
\end{tikzcd}$$
and the canonical extension $\alg{G}_Q' \in \Cat{CExt}_\OO(\alg{G}_{\SC}, \alg{K}_2)$ (defined over $\OO$), exhibit an isomorphism from $\sheaf{D}_Q$ to $\eta^\ast \sheaf{D}_Q^\circ$.  Theorem \ref{ExactBD} gives a trivialization of $\Fun{val}(\sheaf{D}_Q)$ (which depends only on $Q$ restricted to $\sheaf{Y}_{\SC}$),
$$\delta_Q \From \Fun{val}(\sheaf{D}_Q) \xrightarrow{\sim} \sheaf{Y}_{\SC} \oplus \ZZ.$$

Applying this trivialization and the functor $\Fun{val}$ gives an additive functor
$$\Fun{val} \From \Cat{BD}_F(\alg{G}_F, \alg{T}_F) \To \Cat{EZ}_\kk(\alg{\bar G}, \alg{\bar T}),$$
$$(Q, \sheaf{D}, \sheaf{D}_Q \xrightarrow{f} \sheaf{D} ) \mapsto \left( \Fun{val}(\sheaf{D}), \sheaf{Y}_{\SC} \oplus \ZZ \xrightarrow{\Fun{val} f \circ \delta_Q^{-1}} \Fun{val}(\sheaf{D}) \right).$$
\begin{remark}
In order to check that $\Fun{val}$ above is an additive functor, observe that it arises from an additive functor for tori (compatible with pullbacks) and one may check compatibility with Baer sums by checking that the following diagram commutes.
$$\begin{tikzcd}[column sep = 4em]
\Fun{val}(\sheaf{D}_{Q_1}) \Baer \Fun{val}(\sheaf{D}_{Q_2}) \arrow{r}{\delta_{Q_1} \Baer \delta_{Q_2}} \arrow{d}{\sim} & \left( \sheaf{Y}_{\SC} \oplus \ZZ \right) \Baer \left( \sheaf{Y}_{\SC} \oplus \ZZ \right) \arrow{d}{\sim} \\
\Fun{val}(\sheaf{D}_{Q_1 + Q_2}) \arrow{r}{\delta_{Q_1 + Q_2}} & \sheaf{Y}_{\SC} \oplus \ZZ
\end{tikzcd}$$
This can be verified by noting that $\alg{G}_{Q_1}' \Baer \alg{G}_{Q_2}'$ is uniquely isomorphic to $\alg{G}_{Q_1 + Q_2}'$ in the category $\Cat{CExt}_\OO(\alg{G}_{\SC}, \alg{K}_2)$.
\end{remark}

\begin{thm}
\label{BDExactG}
The following sequence of Picard categories and additive functors is exact.
$$\begin{tikzcd}[row sep = 1em]
0 \arrow{r} & \Cat{BD}_\OO(\alg{G}, \alg{T}) \arrow{r}{\eta^\ast} & \Cat{BD}_F(\alg{G}_F, \alg{T}_F) \arrow{r}{\Fun{val}} & \Cat{EZ}_\kk(\alg{\bar G}, \alg{\bar T}) \\
\phantom{0} & (Q, \sheaf{D}^\circ, f^\circ) \arrow[mapsto]{r} & (Q, \eta^\ast \sheaf{D}^\circ, \eta^\ast f^\circ) & \phantom{0} \\
\phantom{0} & \phantom{0} & (Q, \sheaf{D}, f) \arrow[mapsto]{r} & \left( \sheaf{Y}_{\SC} \oplus \ZZ \xrightarrow{\Fun{val} f \circ \delta_Q^{-1}} \Fun{val} ( \sheaf{D} ) \right) 
\end{tikzcd}$$
\end{thm}
\proof
Fix $Q$ for now, which determines an extension $\sheaf{G}_{\mult} \Into \sheaf{D}_Q^\circ \Onto \sheaf{Y}_{\SC}$.  To give a morphism $\sheaf{D}_Q^\circ \xrightarrow{f^\circ} \sheaf{D}^\circ$ of sheaves on $\OO_{\et}$ making the diagram \eqref{IntDDiagram}, it is equivalent by Theorem \ref{ExactBD} (and its compatibility with pullbacks) to give a corresponding morphism of sheaves on $F_{\et}$,
$$\sheaf{D}_Q \xrightarrow{f} \sheaf{D},$$
endowed with an isomorphism $\epsilon$ making the following diagram commute.
$$\begin{tikzcd}
\Fun{val}(\sheaf{D}_Q) \arrow{r}{f} \arrow{d}{\delta_Q} & \Fun{val}(\sheaf{D}) \arrow{d}{\epsilon} \\
\sheaf{Y}_{\SC} \oplus \ZZ \inarrow{r} & \sheaf{Y} \oplus \ZZ
\end{tikzcd}$$
Allowing $Q$ to vary, we find that to give an object of $\Cat{BD}_\OO(\alg{G}, \alg{T})$, it is equivalent to give an object of $\Cat{BD}_F(\alg{G}_F, \alg{T}_F)$ endowed with a trivialization of its image under $\Fun{val}$.  Further details are left to the reader.
\qed

We have now completed a square of Picard categories and additive functors.
\begin{equation}
\tag{Square for $\alg{G}$}
\begin{tikzcd}
\Cat{CExt}_F(\alg{G}_F, \alg{K}_2) \arrow{r}{\partial} \arrow{d}{\Fun{BD}_F} & \Cat{CExt}_{\kk}(\alg{\bar G}, \alg{\bar G}_{\mult}) \arrow{d}{\Fun{EZ}_{\kk}} \\
\Cat{BD}_F(\alg{G}_F, \alg{T}_F) \arrow{r}{\Fun{val}}  & \Cat{EZ}_\kk(\alg{\bar G}, \alg{\bar T})
\end{tikzcd}
\end{equation}

\subsection{The natural isomorphism for tori.}

Let $\alg{T}$ be a split (for now) torus over $\OO$, with character lattice $X$ and cocharacter lattice $Y$.  Consider the square of Picard categories and additive functors.
\begin{equation}
\tag{Square for $\alg{T}$}
\begin{tikzcd}
\Cat{CExt}_F(\alg{T}_F, \alg{K}_2) \arrow{r}{\partial} \arrow{d}{\Fun{BD}_F} & \Cat{Ext}_{\kk}(\alg{\bar T}, \alg{\bar G}_{\mult}) \arrow{d}{\Fun{EZ}_{\kk}} \\
\Cat{BD}_F(\alg{T}_F) \arrow{r}{\Fun{val}} & \Cat{EZ}_\kk(\alg{\bar T})
\end{tikzcd}
\end{equation}
Observe here that $\Cat{CExt}_{\kk}(\alg{\bar T}, \alg{\bar G}_{\mult}) = \Cat{Ext}_{\kk}(\alg{\bar T}, \alg{\bar G}_{\mult})$; every central extension of tori over a field is commutative.

We begin by identifying a natural isomorphism $N \From \Fun{val} \circ \Fun{BD}_F \xRightarrow{\sim} \Fun{EZ}_\kk \circ \partial$, making the square of Picard categories and additive functors 2-commute.  Start with an extension $\alg{K}_2 \Into \alg{T}_C' \Onto \alg{T}$ over $F$, incarnated by an element $C \in X \otimes X$ as in \cite[\S 3.10, 3.11]{B-D}.  Thus, if $C = \sum_{ij} c_{ij} x_i \otimes x_j$, then $\alg{T}_C' = \alg{T} \times \alg{K}_2$ (as a Zariski sheaf) with multiplication given by the rule
$$(s, \alpha) \cdot (t,\beta) = \left( s t, \alpha \beta \cdot \prod_{i,j} \{ x_i(s), x_j(t) \}^{c_{ij}} \right).$$
We trace through the effect of $\Fun{EZ}_\kk \circ \partial$ on the object $\alg{T}_C' \in \Cat{CExt}_F(\alg{T}_F, \alg{K}_2)$.

The pushout of $\alg{T}_C'(F)$ via the tame symbol $\partial \From \alg{K}_2(F) \To \kk^\times$, and restriction to $T^\circ = \alg{T}(\OO)$, gives an extension $\tilde T^\circ = T^\circ \times \kk^\times$, with
$$(u, z_1) \cdot (v, z_2) = \left( uv, z_1 z_2 \cdot \prod_{i,j} \partial \{ x_i(u), x_j(v) \}^{c_{ij}} \right) = (uv, z_1 z_2).$$
In other words, $\tilde T^\circ = T_C^\circ \times \kk^\times$ is a central extension endowed with a splitting.  

From \cite[Construction 12.11]{B-D}, the extension $\tilde T_C^\circ$ is the pullback of the $\kk^\times$-points of the residual extension $\alg{\bar T}_C' = \partial \alg{T}_C'$.  Since we begin with a split torus, and the computations above are compatible with taking \'etale extensions of $\OO$, we find that the residual extension $\alg{\bar T}_C'$ is also equipped with a splitting:
$$\alg{\bar T}_C' = \partial \alg{T}_C' = \alg{\bar T} \times \alg{\bar G}_{\mult}.$$
Applying the functor $\Fun{EZ}_{\kk}$ yields
$$\Fun{EZ}_{\kk} \circ \partial \alg{T}_C' = Y \oplus \ZZ.$$

Next, we trace through the effect of $\Fun{val} \circ \Fun{BD}_F$ on $\alg{T}_C'$.  Its Brylinski-Deligne invariants are $(Q, D_C)$ where $Q(y) = C(y,y)$ (see \cite[\S 3]{B-D}) and $F^\times \Into D_C \Onto Y$ is the extension given by $D_C = Y \times F^\times$ with multiplication
$$(y_1, u_1) \cdot (y_2, u_2) = \left( y_1 + y_2, u_1 u_2 \prod_{i,j} (-1)^{c_{ij} \langle x_i, y_1 \rangle \langle x_j, y_2 \rangle} \right).$$
Pushing out this extension via $\val_\ast$ as in the diagram \eqref{valsplit}, we obtain an extension $\ZZ \Into \val_\ast D_C \Onto Y$ given by $\val_\ast D_C = Y \times \ZZ$ with multiplication
$$(y_1, a_1) \cdot (y_2, a_2) = (y_1 + y_2, a_1 a_2).$$
Indeed, $\val(\pm 1) = 0$, and so $\val_\ast D_C = Y \oplus \ZZ$.

In this way, $\Id \From Y \oplus \ZZ \To Y \oplus \ZZ$ gives an isomorphism
$$N_C \From \Fun{val} \circ \Fun{BD}_F (\alg{T}_C') \xrightarrow{\sim} \Fun{EZ}_\kk \circ \partial (\alg{T}_C').$$

Next, we trace an automorphism of $\alg{T}_C'$ around the diagram.  Following \cite[\S 3.11]{B-D}, every automorphism of $\alg{T}_C'$ arises from an element of $X \otimes F^\times$.  For $x \otimes s \in X \otimes F^\times$, the corresponding automorphism $\alpha_{x \otimes s}$ of $\alg{T}_C'$ is given by
$$\alpha_{x \otimes s}(t, \kappa) = (t, \kappa \cdot \{ x(t), s \} ).$$
We trace this automorphism through the functor $\Fun{EZ}_\kk \circ \partial$ first.  On the (split) extension $\tilde T_C^\circ = T^\circ \times \kk^\times$, we obtain the automorphism
$$\alpha_{x \otimes s}(u, z) = (u,z \cdot \partial \{x(u), s \} ) = \left( u,z \cdot \overline{x(u)}^{\val(s)} \right).$$
Here $u \in T^\circ$, and $x(u) \in \OO^\times$ has reduction $\overline{x(u)} \in \kk^\times$.

Passing to the residual extension gives an automorphism of $\alg{\bar T}_C' = \alg{\bar T} \times \alg{\bar G}_{\mult}$,
$$\alpha_{x \otimes s}(\bar u, z) = \left( \bar u, z \cdot x(\bar u)^{\val(s)} \right), \text{ for all } \bar u \in \alg{\bar T}, z \in \alg{\bar G}_{\mult}.$$
On the cocharacter lattice $Y \oplus \ZZ$ of $\alg{\bar T} \times \alg{\bar G}_{\mult}$, we find the automorphism
\begin{equation}
\label{EAut}
\alpha_{x \otimes s}(y, a) = (y, a + \val(s) \langle x,y \rangle).
\end{equation}

We trace the automorphism of $\alg{T}_C'$ through the functor $\Fun{val} \circ \Fun{BD}_F$ now.  On the Brylinski-Deligne invariant $F^\times \Into D_C \Onto Y$, we find that
$$\alpha_{x \otimes s}(y, u) = (y, u \cdot s^{\langle x, y \rangle} )$$
following the remarks of \cite[\S 3.11]{B-D}.  Pushing out yields an automorphism of $\ZZ \Into \val_\ast D_C \Onto Y$,
\begin{equation}
\label{BAut}
\alpha_{x \otimes s}(y, a) = (y, a + \langle x,y \rangle \val(s) ).
\end{equation}

The coincidences between \eqref{EAut} and \eqref{BAut} yield the following theorem.
\begin{thm}
\label{ResExtSplitTori}
For every split algebraic torus $\alg{T}$ over $\OO$, with cocharacter lattice $Y$, there exists a unique natural isomorphism of additive functors $N_Y \From \Fun{val} \circ \Fun{BD}_F \xRightarrow{\sim} \Fun{EZ} \circ \partial$, agreeing with $N_C$ on the extensions $\alg{T}_C'$, making the square below 2-commute.
\begin{equation}
\tag{Square for $\alg{T}$}
\begin{tikzcd}[row sep=3em]
\Cat{CExt}_F(\alg{T}_F, \alg{K}_2) \arrow{r}{\partial} \arrow{d}{\Fun{BD}_F} & \Cat{Ext}_{\kk}(\alg{\bar T}, \alg{\bar G}_{\mult}) \arrow{d}{\Fun{EZ}_\kk} \\
\Cat{BD}_F(\alg{T}_F) \arrow{r}{\Fun{val}} \arrow[Rightarrow, bend right=10]{ur}{N_Y} & \Cat{EZ}_\kk(\alg{\bar T})
\end{tikzcd}
\end{equation}
\end{thm}
\proof
By \cite[Proposition 3.11]{B-D}, every object of $\Cat{CExt}_F(\alg{T}_F, \alg{K}_2)$ is isomorphic to $\alg{T}_C'$ for some $C \in X \otimes X$.  For such an object $\alg{T}'$, choose an isomorphism $\iota \From \alg{T}' \To \alg{T}_C'$ for some $C \in X \otimes X$.  Then there is a unique isomorphism $N_\iota$ making the following diagram in $\Cat{EZ}_\kk(\alg{\bar T}) = \Cat{Ext}(Y, \ZZ)$ commute.
$$\begin{tikzcd}
\Fun{val} \circ \Fun{BD}_F(\alg{T}') \arrow{r}{N_\iota}[swap]{\sim} \arrow{d}{\sim}[swap]{\Fun{val} \circ \Fun{BD}_F(\iota)} & \Fun{EZ}_\kk \circ \partial(\alg{T}') \arrow{d}{\Fun{EZ}_\kk \circ \partial(\iota)}[swap]{\sim} \\
\Fun{val} \circ \Fun{BD}_F(\alg{T}_C') \arrow{r}{N_C}[swap]{\sim} & \Fun{EZ}_\kk \circ \partial(\alg{T}_C')
\end{tikzcd}$$

If another isomorphism $\jmath \From \alg{T}' \To \alg{T}_C'$ is chosen, one obtains another isomorphism $N_\jmath$ accordingly.  But then $\jmath = \alpha \circ \iota$ for some $\alpha \in \Aut(\alg{T}_C') = X \otimes F^\times$.  The coincidence between \eqref{EAut} and \eqref{BAut} implies that $N_{\jmath} = N_\iota$.

Now consider another element $C_0 \in X \otimes X$ such that $\alg{T}_{C_0}$ is isomorphic to $\alg{T}'$ as well.  Define the bilinear form $A = C - C_0 \in X \otimes X$.  Since $\alg{T}_C'$ is isomorphic to $\alg{T}_{C_0}'$, the associated quadratic forms $y \mapsto C(y,y)$ and $y \mapsto C_0(y,y)$ are equal; thus $A(y,y) = 0$ for all $y \in Y$.  

If $A = \sum_{i,j} a_{ij} x_i \otimes x_j$, then an isomorphism $\alpha \From \alg{T}_{C_0}' \xrightarrow{\sim} \alg{T}_C'$ is given by the formula
$$\alpha(t, \kappa) \mapsto \left( t, \kappa \cdot \prod_{i,j} \{ x_i(t), x_j(t) \}^{a_{ij}} \right).$$
Define $\iota_0 = \alpha^{-1} \circ \iota \From \alg{T}' \To \alg{T}_{C_0}'$.  We have the following diagram in $\Cat{EZ}_\kk(\alg{\bar T})$.
$$\begin{tikzcd}[row sep=3em, column sep = 4.5em]
\Fun{val} \circ \Fun{BD}_F(\alg{T}_{C_0}') \arrow{r}{N_{C_0}}[swap]{\sim} \arrow[bend right=70]{dd}[swap]{\Fun{val} \circ \Fun{BD}_F(\alpha)} & \Fun{EZ}_\kk \circ \partial(\alg{T}_{C_0}') \arrow[bend left=70]{dd}{\Fun{EZ}_\kk \circ \partial(\alpha)} \\
\Fun{val} \circ \Fun{BD}_F(\alg{T}') \arrow{r}{N_\iota = N_{\iota_0}?}[swap]{\sim} \arrow{d}{\Fun{val} \circ \Fun{BD}_F(\iota)}[swap]{\sim} \arrow{u}{\sim}[swap]{\Fun{val} \circ \Fun{BD}_F(\iota_0)} & \Fun{EZ}_\kk \circ \partial(\alg{T}') \arrow{d}{\sim}[swap]{\Fun{EZ}_\kk \circ \partial(\iota)} \arrow{u}{\Fun{EZ}_\kk \circ \partial(\iota_0)}[swap]{\sim}  \\
\Fun{val} \circ \Fun{BD}_F(\alg{T}_C') \arrow{r}{N_C}[swap]{\sim} & \Fun{EZ}_\kk \circ \partial(\alg{T}_C')
\end{tikzcd}$$
The top square commutes when choosing the morphism $N_{\iota_0}$ (along the middle row) and the bottom square commutes when choosing the morphism $N_\iota$.

Now we trace through $\Fun{EZ}_\kk \circ \partial(\alpha)$, an isomorphism from $Y \oplus \ZZ$ to itself.  First we obtain an isomorphism from $\tilde T_C^\circ = T^\circ \times \kk^\times$ to $\tilde T_{C_0}^\circ = T^\circ \times \kk^\times$, given by
$$\alpha(u,z) = \left( u, z \cdot \prod_{i,j} \partial \{ x_i(u), x_j(u) \}^{a_{ij}} \right) = (u,z).$$
Thus the isomorphism $\Fun{EZ}_\kk \circ \partial(\alpha)$ is the identity on $Y \oplus \ZZ$.  A short computation demonstrates that $\Fun{val} \circ \Fun{BD}_F(\alpha)$ is the identity on $Y \oplus \ZZ$ as well.

It follows that $\Fun{val} \circ \Fun{BD}_F(\iota_0) = \Fun{val} \circ \Fun{BD}_F(\iota)$, and similarly $\Fun{EZ}_\kk \circ \partial(\iota_0) = \Fun{EZ}_\kk \circ \partial(\iota)$; thus the isomorphism $N_{\iota_0}$ which makes the top square commute coincides with the isomorphism $N_\iota$ which makes the bottom square commute.

Thus we find an isomorphism
$$N_Y[\alg{T}'] \From \Fun{val} \circ \Fun{BD}_F(\alg{T}') \xrightarrow{\sim} \Fun{EZ}_\kk \circ \partial(\alg{T}'),$$
defined independently of the choice of $C \in X \otimes X$ and the choice of isomorphism $\alg{T}' \To \alg{T}_C'$.  This suffices to define a unique natural isomorphism of functors $N_Y \From \Fun{val} \circ \Fun{BD}_F \xRightarrow{\sim} \Fun{EZ} \circ \partial$.  Naturality follows from compatibility with automorphisms of each $\alg{T}_C'$.  The natural automorphism respects the additive structure of the Picard categories, since the ``sum'' of extensions corresponds to addition in $X \otimes X$, and one may check compatibility for incarnated extensions such as $\alg{T}_C'$.  
\qed

Suppose now that $\tau \From \alg{T}_1 \To \alg{T}_2$ is an isomorphism of split tori, corresponding to a isomorphism $\tau \From Y_1 \To Y_2$ of cocharacter lattices and $\tau^\ast \From X_2 \To X_1$ of character lattices.  The isomorphism $\tau$ defines pullback functors:
\begin{align*}
\tau^\ast & \From \Cat{CExt}(\alg{T}_{2,F}, \alg{K}_2) \To \Cat{CExt}(\alg{T}_{1,F}, \alg{K}_2); \\
\tau^\ast & \From \Cat{Ext}_\kk(\alg{\bar T}_2, \alg{\bar G}_{\mult}) \To \Cat{Ext}_\kk(\alg{\bar T}_1, \alg{\bar G}_{\mult}); \\ 
\tau^\ast & \From \Cat{BD}_F(\alg{T}_{2,F}) \To  \Cat{BD}_F(\alg{T}_{1,F}); \\
\tau^\ast & \From \Cat{EZ}_\kk(\alg{\bar T}_2) \To \Cat{EZ}_\kk(\alg{\bar T}_1).
\end{align*}
The inverse $\tau^{-1}$ defines pullback functors in the other direction, with natural isomorphisms $(\tau^{-1})^\ast \tau^\ast \Rightarrow \Id$.  

We find a cube of Picard categories and additive functors.
$$\begin{tikzcd}[back line/.style={densely dotted}, row sep=3em, column sep=1.5em]
& \Cat{CExt}_F(\alg{T}_{2,F}, \alg{K}_2)  \ar{dl}[swap, sloped]{\tau^\ast} \ar{rr}{\partial_2} \ar[back line]{dd}[near end]{\Fun{BD}_{2,F}} 
  & & \Cat{Ext}_{\kk}(\alg{\bar T}_2, \alg{\bar G}_{\mult}) \ar{dd}{\Fun{EZ}_{2,\kk}} \ar{dl}[swap,sloped,near start]{\tau^\ast} \\
\Cat{CExt}_F(\alg{T}_{1,F}, \alg{K}_2) \ar[crossing over]{rr}[near start]{\partial_1} \ar{dd}[swap]{\Fun{BD}_{1,F}} 
  & & \Cat{Ext}_{\kk}(\alg{\bar T}_{1}, \alg{\bar G}_{\mult}) \\
& \Cat{BD}_F(\alg{T}_{2,F}) \ar[back line]{rr}[near start]{\Fun{val}_2} \ar[back line]{dl}[swap,sloped,near start]{\tau^\ast}
  & & \Cat{EZ}_\kk(\alg{\bar T}_2) \ar{dl}[swap,sloped,near start]{\tau^\ast} \\
\Cat{BD}_F(\alg{T}_{1,F}) \ar{rr}{\Fun{val}_1} & &  \Cat{EZ}_\kk(\alg{\bar T}_{1}) \ar[crossing over, leftarrow]{uu}[swap, near start]{\Fun{EZ}_{1,\kk}}
\end{tikzcd}
$$
Pullbacks along $\tau$ commute with the functors $\partial$, $\Fun{EZ}_\kk$, $\Fun{BD}_F$, and $\Fun{val}$, up to natural isomorphism.  In other words, there are natural isomorphisms $\tau^\ast \partial_2 \Leftrightarrow \partial_1 \tau^\ast$, etc..  Thus the lateral faces of the cube commute up to natural isomorphism.   For ease of notation, define $\eta_i = \Fun{EZ}_{i, \kk} \circ \partial_i$ and $\beta_i = \Fun{val}_i \circ \Fun{BD}_{i,F}$ for $i = 1,2$.  Compatibility with pullbacks gives natural isomorphisms,
$$I_\beta \From \beta_1 \Rightarrow \tau^\ast \beta_2 (\tau^{-1})^\ast, \quad I_\eta \From \eta_1 \Rightarrow \tau^\ast \eta_2 (\tau^{-1})^\ast.$$

The front face of the cube 2-commutes via the natural isomorphism $N_{Y_1}$ and the back face 2-commutes via $N_{Y_2}$, according to the previous theorem.  Note that $N_{Y_i} \From \beta_i \Rightarrow \eta_i$ is a natural isomorphism in our abbreviated notation.  These natural isomorphisms are compatible with pullbacks in the following sense.
\begin{proposition}
The following diagram commutes.
$$\begin{tikzcd}
\beta_1 \arrow{r}{N_{Y_1}} \arrow{d}{I_\beta} & \eta_1 \arrow{d}{I_\eta} \\
\tau^\ast \beta_2 (\tau^{-1})^\ast \arrow{r}{N_{Y_2}} & \tau^\ast \eta_2 (\tau^{-1})^\ast
\end{tikzcd}$$
Here the entries in the diagram are objects and morphisms in the Picard category of additive functors, $\Fun{Hom} \left( \Cat{CExt}_F(\alg{T}_{1,F}, \alg{K}_2), \Cat{EZ}_\kk(\alg{\bar T}_1) \right)$.
\end{proposition}
\proof
Since the natural isomorphisms are uniquely determined by their behavior on ``incarnated'' extensions, it suffices to consider $C \in X_2 \otimes X_2$ with $C_1 = \tau^\ast(C_2) \in X_1 \otimes X_1$; we have corresponding objects $\alg{T}_{C_2}' \in \Cat{CExt}_F(\alg{T}_{2,F}, \alg{K}_2)$ and $\alg{T}_{C_1}' \in \Cat{CExt}_F(\alg{T}_{1,F}, \alg{K}_2)$.  A direct computation demonstrates that
$$\tau^\ast \alg{T}_{C_2}' = \alg{T}_{C_1}', \quad (\tau^{-1})^\ast \alg{T}_{C_1}' = \alg{T}_{C_2}'.$$

We find that $\beta_1(\alg{T}_{C_1}') = Y_1 \oplus \ZZ$, and similarly,
$$\tau^\ast \beta_2 (\tau^{-1})^\ast \alg{T}_{C_1}' = \tau^\ast \beta_2(\alg{T}_{C_2}') = \tau^\ast(Y_2 \oplus \ZZ) = Y_1 \oplus \ZZ.$$
Thus the natural isomorphism $I_\beta$ is given by the identity map from $Y_1 \oplus \ZZ$ to itself.  Similarly, $I_\eta$ is given by the identity map from $Y_1 \oplus \ZZ$ to itself.  The natural isomorphism $N_{Y_2}$ is the identity map on $Y_2 \oplus \ZZ$, and the natural isomorphism $N_{Y_1}$ is the identity on $Y_1 \oplus \ZZ$.  
It follows quickly that the diagram commutes.
\qed

For nonsplit tori, we may exploit the fact that all the Picard categories discussed here arise from Picard {\em stacks} on $\OO_{\et}$.  Consider a torus $\alg{T}$ over $\OO$, not necessarily split.  For a finite \'etale $\OO' / \OO$, with $\Spec(\OO')$ connected, write $F' / F$ for the resulting unramified extension of fraction fields, and $\kk' / \kk$ for the extension of residue fields.  For a Picard stack $\stack{P}$ over $\OO_{\et}$, write $\stack{P}[\OO']$ for the Picard category of sections of $\stack{P}$ over $\OO'$.  

Define four Picard stacks on $\OO_{\et}$ as follows:
\begin{align*}
\stack{CExt}(\alg{T}_F, \alg{K}_2)[\OO'] &\defeq \Cat{CExt}_{F'}(\alg{T}_{F'}, \alg{K}_2); \\
\stack{Ext}(\alg{\bar T}, \alg{\bar G}_{\mult})[\OO'] &\defeq \Cat{Ext}_{\kk'}(\alg{\bar T}_{\kk'}, \alg{\bar G}_{\mult}); \\
\stack{BD}(\alg{T}_F)[\OO'] &\defeq \Cat{BD}_{F'}(\alg{T}_{F'}); \\
\stack{EZ}(\alg{\bar T})[\OO'] &\defeq \Cat{EZ}_{\kk'}(\alg{\bar T}_{\kk'}).
\end{align*}
The base extension functors for these stacks have already been discussed.

When $\stack{P}_1$ and $\stack{P}_2$ are two Picard stacks on $\OO_{\et}$, the additive functors from $\stack{P}_1$ to $\stack{P}_2$ and natural transformations between them form a Picard stack $\stack{Hom}(\stack{P}_1, \stack{P}_2)$.  Thus, as the functors $\partial$, $\Fun{BD}$, $\Fun{EZ}$, $\Fun{val}$, and the natural isomorphisms $N$, are all defined in such a way to be compatible with base extensions $\OO' / \OO$ and automorphisms of split tori, descent gives a 2-commutative square of Picard stacks on $\OO_{\et}$.
\begin{equation}
\tag{Square for $\alg{T}$}
\begin{tikzcd}[row sep=3em]
\stack{CExt}(\alg{T}_F, \alg{K}_2) \arrow{r}{\partial} \arrow{d}{\Fun{BD}} & \stack{Ext}(\alg{\bar T}, \alg{\bar G}_{\mult}) \arrow{d}{\Fun{EZ}} \\
\stack{BD}(\alg{T}_F) \arrow{r}{\Fun{val}} \arrow[Rightarrow, bend right=10]{ur}{N_Y} & \stack{EZ}(\alg{\bar T})
\end{tikzcd}
\end{equation}

In this way, descent extends Theorem \ref{ResExtSplitTori} to the nonsplit case.
\begin{corollary}
For every torus $\alg{T}$ over $\OO$, with cocharacter lattice $\sheaf{Y}$, descent provides a natural isomorphism of functors $N_{\sheaf{Y}} \From \Fun{val} \circ \Fun{BD_F} \To \Fun{EZ}_\kk \circ \partial$ making the square of Picard categories and additive functors 2-commute.
\begin{equation}
\tag{Square for $\alg{T}$}
\begin{tikzcd}[row sep=3em]
\Cat{CExt}_F(\alg{T}_F, \alg{K}_2) \arrow{r}{\partial} \arrow{d}{\Fun{BD}_F} & \Cat{Ext}_{\kk}(\alg{\bar T}, \alg{\bar G}_{\mult}) \arrow{d}{\Fun{EZ}_\kk} \\
\Cat{BD}_F(\alg{T}_F) \arrow{r}{\Fun{val}} \arrow[Rightarrow, bend right=10]{ur}{N_{\sheaf{Y}}} & \Cat{EZ}_\kk(\alg{\bar T})
\end{tikzcd}
\end{equation}
These natural isomorphisms are compatible with pullbacks for isomorphisms $\alg{T}_1 \To \alg{T}_2$ of tori over $\OO$.
\end{corollary}

\subsection{Reductive groups}

Now we consider the general case of a reductive group $\alg{G}$ over $\OO$ endowed with a maximal torus $\alg{T}$ over $\OO$.  As before $\jota \From \alg{G}_{\SC} \To \alg{G}$ denotes the simply-connected cover of the derived subgroup.  We have constructed a square of Picard categories and additive functors.
\begin{equation}
\tag{Square for $\alg{G}$}
\begin{tikzcd}
\Cat{CExt}_F(\alg{G}_F, \alg{K}_2) \arrow{r}{\partial} \arrow{d}{\Fun{BD}_F} & \Cat{CExt}_{\kk}(\alg{\bar G}, \alg{\bar G}_{\mult}) \arrow{d}{\Fun{EZ}_{\kk}} \\
\Cat{BD}_F(\alg{G}_F, \alg{T}_F) \arrow{r}{\Fun{val}} & \Cat{EZ}_\kk(\alg{\bar G}, \alg{\bar T})
\end{tikzcd}
\end{equation}

Define $\beta = \Fun{val} \circ \Fun{BD}_F$ and $\eta = \Fun{EZ}_\kk \circ \partial$.  In the case of tori and simply-connected semisimple groups, we have constructed a natural isomorphism $N \From \beta \To \eta$.  Here we construct such a natural isomorphism $N$ in general.

Begin with an object $\alg{G}_F' \in \Cat{CExt}_F(\alg{G}_F, \alg{K}_2)$.  Write $\alg{\bar G}' = \partial \alg{G}_F'$ for its residual extension.  Write $(Q, \sheaf{D}, f)$ for its Brylinski-Deligne invariants.  

Write $\alg{T}_F' \in \Cat{CExt}_F(\alg{G}_F, \alg{K}_2)$ for the pullback of $\alg{G}_F'$, and $\alg{\bar T}'$ for its residual extension.  Thus $(Q, \sheaf{D}) = \Fun{BD}_F(\alg{T}_F')$.  

Let $\sheaf{Y}' \in \Cat{Ext}_{\kk_{\et}}(\sheaf{Y}, \ZZ)$ denote the cocharacter lattice of the torus $\alg{\bar T}'$, i.e., $\sheaf{Y}' = \Fun{EZ}_\kk(\alg{\bar T}')$.  Observe that $\Fun{val}(\sheaf{D}) = \beta \alg{T}_F'$ and $\sheaf{Y}' = \eta \alg{T}_F'$.

Write $\sheaf{Y}_{\SC}'$ for the pullback of $\sheaf{Y}'$ to $\sheaf{Y}_{\SC}$.  The unique splitting $\alg{\bar G}_{\SC}' \isom \alg{\bar G}_{\SC} \times \alg{\bar G}_{\mult}$ provides an isomorphism
$$\can \From \sheaf{Y}_{\SC} \oplus \ZZ \xrightarrow{\sim} \sheaf{Y}_{\SC}'.$$
To give an isomorphism, $N_{\alg{G}} \From \beta \alg{G}_F' \To \eta \alg{G}_F'$, is to give a isomorphism $\nu \From \Fun{val}(\sheaf{D}) \To \sheaf{Y}'$ in $\Cat{CExt}_{\kk_{\et}}(\sheaf{Y}, \ZZ)$, making the diagram commute.
\begin{equation}
\label{sstCompat}
\begin{tikzcd}[column sep = 3.5em]
\sheaf{Y}_{\SC} \oplus \ZZ \arrow{r}{\Fun{val} f \circ \delta_Q^{-1}}[swap]{\sim} \arrow{d}{=}  & \Fun{val} (\sheaf{D}_{\SC}) \inarrow{r} & \Fun{val}(\sheaf{D}) \arrow{d}{\nu} \\
\sheaf{Y}_{\SC} \oplus \ZZ \arrow{r}{\can}[swap]{\sim} & \sheaf{Y}_{\SC}' \inarrow{r} & \sheaf{Y}'
\end{tikzcd}
\end{equation}

We have already constructed an isomorphism $N_{\sheaf{Y}} \From \Fun{val}(\sheaf{D}) \To \sheaf{Y}'$; the reductive group case rests now on compatibility of this toral case with the simply-connected semisimple case.
\begin{lemma}
If $\nu = N_{\sheaf{Y}}$, then Diagram \eqref{sstCompat} commutes.
\end{lemma} 
\proof
Restricting $N_{\sheaf{Y}}$ to $\Fun{val}(\sheaf{D}_{\SC})$ gives an isomorphism in $\Cat{Ext}_{\kk_{\et}}(\sheaf{Y}_{\SC}, \ZZ)$.  There exists a unique homomorphism $\phi \in \Hom(\sheaf{Y}_{\SC}, \ZZ)$ which makes the following diagram commute.
$$\begin{tikzcd}[column sep = 3.5em]
\sheaf{Y}_{\SC} \oplus \ZZ \arrow{r}{\Fun{val} f \circ \delta_Q^{-1}}[swap]{\sim} \arrow{d}[swap]{(y,n) \mapsto (y,\phi(y) + n)}  & \Fun{val} (\sheaf{D}_{\SC}) \inarrow{r} \arrow{d}{N_{\sheaf{Y}}}[swap]{\sim} & \Fun{val}(\sheaf{D}) \arrow{d}{N_{\sheaf{Y}}}[swap]{\sim} \\
\sheaf{Y}_{\SC} \oplus \ZZ \arrow{r}{\can}[swap]{\sim} & \sheaf{Y}_{\SC}' \inarrow{r} & \sheaf{Y}'
\end{tikzcd}$$

We claim that $\phi = 0$; to prove this it suffices to work \'etale locally, and assume that $\sheaf{Y}$ and the Weyl group $\sheaf{W}$ are constant sheaves.  Let $\dot w$ be an element of $\alg{G}_{\SC}(\OO)$ lying over an element $w \in \sheaf{W}[\OO]$.  Conjugation by $\dot w$ gives a homomorphism of $\OO$-tori, $\Int(w) \From \alg{T} \To \alg{T}$, and a corresponding homomorphism $\Int(w) \From \sheaf{Y} \To \sheaf{Y}$.

We may use $\Int(w)$ to pull back central extensions and define
$$\Int(w)^\ast \alg{T}_F' \in \Cat{CExt}(\alg{T}_F, \alg{K}_2), \quad \Int(w)^\ast \sheaf{D} \in \Cat{CExt}_{F_{\et}}(\sheaf{Y}, \sheaf{G}_{\mult}), \quad \text{etc..}$$

On the other hand, the representative $\dot w \in \alg{G}_{\SC}(\OO)$ gives a (lifting-then) conjugation map $\Int(\dot w)$ fitting into a commutative diagram
$$\begin{tikzcd}
\alg{K}_2 \inarrow{r} \arrow{d}{=} & \alg{G}' \onarrow{r} \arrow{d}{\Int(\dot w)} & \alg{G} \arrow{d}{\Int(\dot w)} \\ 
\alg{K}_2 \inarrow{r} & \alg{G}' \onarrow{r} & \alg{G} 
\end{tikzcd}$$
Restricting this to $\alg{T}$, we find that $\Int(\dot w)$ gives an isomorphism,
$$\Int(\dot w) \From \alg{T}_F' \xrightarrow{\sim} \Int(w)^\ast \alg{T}_F', \quad \text{ in } \Cat{CExt}_F(\alg{T}_F, \alg{K}_2).$$

As the natural isomorphism $N_{\sheaf{Y}}$ and functors $\eta, \beta$, are compatible with pullbacks, we find a commutative square
$$\begin{tikzcd}
\Fun{val}(\sheaf{D}) \arrow{r}{\Int(\dot w)} \arrow{d}{N_{\sheaf{Y}}} & \Int(w)^\ast \Fun{val}(\sheaf{D}) \arrow{d}{\Int(w)^\ast N_{\sheaf{Y}}} \\
\sheaf{Y}' \arrow{r}{\Int(\dot w)} & \Int(w)^\ast \sheaf{Y}'
\end{tikzcd}$$

Tracing these maps to $\sheaf{Y}_{\SC} \oplus \ZZ$, we find a commutative diagram,
$$\begin{tikzcd}[column sep = 6em]
\sheaf{Y}_{\SC} \oplus \ZZ \arrow{r}{(y,n) \mapsto ({}^w y,n)} \arrow{d}[swap]{(y,n) \mapsto (y,\phi(y) + n)} & \sheaf{Y}_{\SC} \oplus \ZZ \arrow{d}{(y,n) \mapsto (y,\phi(y) + n)} \\
\sheaf{Y}_{\SC} \oplus \ZZ \arrow{r}{(y,n) \mapsto ({}^w y,n)} & \sheaf{Y}_{\SC} \oplus \ZZ
\end{tikzcd}$$
In other words, for all $n \in \ZZ$, $y \in \sheaf{Y}_{\SC}$, and all Weyl elements $w$, we have
$$({}^w y,\phi({}^w y) + n) = ({}^w y, \phi(y) + n).$$
Thus $\phi \From \sheaf{Y}_{\SC} \To \ZZ$ is Weyl-invariant, and therefore $\phi = 0$.
\qed

\begin{thm}
\label{natisomG}
The natural isomorphism $N_{\sheaf{Y}}$, defined earlier for tori, defines a natural isomorphism $N_{\alg{G}}$ making the following diagram of Picard categories and additive functors 2-commute. 
\begin{equation}
\tag{Square for $\alg{G}$}
\begin{tikzcd}
\Cat{CExt}_F(\alg{G}_F, \alg{K}_2) \arrow{r}{\partial} \arrow{d}{\Fun{BD}_F} & \Cat{CExt}_{\kk}(\alg{\bar G}, \alg{\bar G}_{\mult}) \arrow{d}{\Fun{EZ}_{\kk}} \\
\Cat{BD}_F(\alg{G}_F, \alg{T}_F) \arrow{r}{\Fun{val}} \arrow[Rightarrow, bend right=10]{ur}{N_{\alg{G}}} & \Cat{EZ}_\kk(\alg{\bar G}, \alg{\bar T})
\end{tikzcd}
\end{equation}
\end{thm}
\proof
The lemma demonstrates that $N_{\sheaf{Y}}$ provides an isomorphism,
$$N[\alg{G}_F'] \From \beta \alg{G}_F' \xrightarrow{\sim} \eta \alg{G}_F',$$
for all $\alg{G}_F' \in \Cat{CExt}_F(\alg{G}_F, \alg{K}_2)$.  

Next, consider a morphism $c \From \alg{G}_{1,F}' \To \alg{G}_{2,F'}$ in $\Cat{CExt}_F(\alg{G}_F, \alg{K}_2)$.  This yields a diagram in $\Cat{EZ}_{\kk}(\alg{\bar G}, \alg{\bar T})$.
\begin{equation}
\label{natisomN}
\begin{tikzcd}[column sep = 4.5em]
\beta \alg{G}_{1,F}' \arrow{r}{N[\alg{G}_{1,F}']} \arrow{d}{\beta(c)} & \eta \alg{G}_{1,F}' \arrow{d}{\eta(c)} \\
\beta \alg{G}_{2,F}' \arrow{r}{N[\alg{G}_{2,F}']} & \eta \alg{G}_{2,F}' 
\end{tikzcd}
\end{equation}
Pulling back from $\alg{\bar G}$ to $\alg{\bar T}$ gives a faithful additive functor from $\Cat{EZ}_{\kk}(\alg{\bar G}, \alg{\bar T})$ to $\Cat{EZ}_{\kk}(\alg{\bar T}) = \Cat{Ext}(\sheaf{Y}, \ZZ)$.  Thus to check that Diagram \eqref{natisomN} commutes, it suffices to check that it commutes after pulling back to $\alg{T}$ throughout.  But this commutativity follows from the fact that $N_{\sheaf{Y}}$ is a natural isomorphism of functors.  

Hence Diagram \eqref{natisomN} commutes, and we find a natural isomorphism of functors $N_{\alg{G}}$ as desired.  The compatibility with the additive structure follows as well from the case of tori.
\qed

\subsection{The classification theorem}

Here we keep $\alg{G}$, a reductive group over $\OO$ with maximal torus $\alg{T}$ over $\OO$.  The previous sections provide a diagram of Picard categories and additive functors, with exact rows, and natural isomorphisms expressing its 2-commutativity.
$$\begin{tikzcd}
0 \arrow{r} & \Cat{CExt}_{\OO}(\alg{G}, \alg{K}_2) \arrow{r} \arrow{d}{\Fun{BD}_\OO} & \Cat{CExt}_F(\alg{G}_F, \alg{K}_2) \arrow{r}{\partial} \arrow{d}{\Fun{BD}_F} & \Cat{CExt}_{\kk}(\alg{\bar G}, \alg{\bar G}_{\mult}) \arrow{d}{\Fun{EZ}_\kk} \\
0 \arrow{r} & \Cat{BD}_{\OO}(\alg{G}, \alg{T}) \arrow{r} & \Cat{BD}_F(\alg{G}_F, \alg{T}_F) \arrow{r}{\Fun{val}} & \Cat{EZ}(\alg{\bar G}, \alg{\bar T})
\end{tikzcd}$$

Exactness of the top row is a special case of the exact sequence \eqref{CExtExactBetter}.
Exactness of the bottom row is Theorem \ref{BDExactG}.
The 2-commutativity of the right square is Theorem \ref{natisomG}.
The 2-commutativity of the left square is compatibility with pullback, see Diagram \eqref{BDPull}.
The functor $\Fun{BD}_F$ is an equivalence by \cite[Theorem 7.2]{B-D}.
The functor $\Fun{EZ}_\kk$ is an equivalence by Theorem \ref{EZMain}.

It follows that $\Fun{BD}_\OO$ is an equivalence, and so we find a classification which extends the main result of \cite{B-D}.
\begin{thm}
When $\OO$ is a discrete valuation ring, with finite residue field or containing a field, $\Fun{BD}_\OO$ is an equivalence of Picard categories,
$$\Fun{BD}_\OO \From \Cat{CExt}_{\OO}(\alg{G}, \alg{K}_2) \xrightarrow{\sim} \Cat{BD}_{\OO}(\alg{G}, \alg{T}).$$
\end{thm}

\bibliography{CovLang2014.bib}

%\printbibliography

\end{document}